\documentclass[12pt]{article}%
\usepackage{amsfonts}
\usepackage{fancyhdr}
\usepackage[colorlinks = true,
            linkcolor = blue,
            urlcolor  = blue,
            citecolor = black,
            anchorcolor = blue]{hyperref}
\usepackage[a4paper, top=2.5cm, bottom=2.5cm, left=2cm, right=2cm]%
{geometry}
\usepackage{indentfirst}
\usepackage{times}
\usepackage{amsmath}
\usepackage{mathtools}
\usepackage{changepage}
\usepackage{enumerate}
\usepackage{tcolorbox}
\usepackage{float}
\usepackage{amssymb}
\usepackage{bbm}
\usepackage{graphicx}%
\usepackage{cleveref}
\newcommand{\rev}{{\text{rev}}}
\newcommand{\starplus}{{\prescript{+}{}{\mathsf{Star}}}}
\newcommand{\starreg}{{\mathsf{Star}}}
\newcommand{\FS}{{\mathsf{FS}}}

\newenvironment{proof}[1][Proof]{\textbf{#1.} }{\ \rule{0.5em}{0.5em}}

\title{On the Asymmetric Generalizations of Two Extremal Questions on Friends-and-Strangers Graphs}
\author{Kiril Bangachev}
\date{}

\begin{document}
\maketitle

\begin{abstract}
\noindent
    For two graphs $X$ and $Y$ with vertex sets $V(X)$ and $V(Y)$ of the same cardinality $n,$ the friends-and-strangers graph $\FS(X,Y)$ was recently defined by Defant and Kravitz. The vertices of $\FS(X,Y)$ are the bijections from $V(X)$ to $V(Y),$ and two bijections $\sigma$ and $\tau$ are adjacent if they agree everywhere except at two vertices 
    $a,b\in V(X)$ such that $a$ and $b$ are adjacent in $X$ and $\sigma(a)$ and $\sigma(b)$ are adjacent in $Y.$ We study generalized versions of two problems by Alon, Defant, and Kravitz. First, we show that if $X$ and $Y$ have minimum degrees $\delta(X)$ and $\delta(Y)$ that satisfy $\delta(X)> n/2, \delta(Y)>n/2,$ and $2\min(\delta(X), \delta(Y))+3\max(\delta(X), \delta(Y))\ge 3n,$ then $\FS(X,Y)$ is connected. As a corollary, we settle a recent conjecture by Alon, Defant, and Kravitz stating that there exists a number $d_n = 3n/5 + O(1)$ such that if both $X$ and $Y$ have minimum degrees at least $d_n,$ the graph $\FS(X,Y)$ is connected. When $X$ and $Y$ are bipartite, a parity obstruction prevents $\FS(X,Y)$ from being connected. We show that 
    if $X$ and $Y$ are edge-subgraphs of $K_{r,r}$ that satisfy $\delta(X)+\delta(Y)\ge 3r/2+1,$ 
    then the graph $\FS(X,Y)$ has exactly two connected components. As a corollary, we provide an almost complete answer to another recent question of Alon, Defant, and Kravitz asking for the minimum number $d^*_{r,r}$ such that for any
    edge-subgraph $X$ of $K_{r,r}$ satisfying $\delta(X)\ge d^*_{r,r},$ the graph $\FS(X,K_{r,r})$ has exactly two connected components. We show that 
    $d^*_{r,r} = r/2+1$ when $r$ is even and $d^*_{r,r}\in \{\lceil r/2\rceil, \lceil r/2\rceil+1\}$ when $r$ is odd.
\end{abstract}

\section{Introduction}
\subsection{Background}
Defant and Kravitz \cite{FSFixed} recently introduced the concept of friends-and-strangers graphs. The formal definition is given as follows: \\

\noindent
\textbf{Definition 1.1} (\cite{FSFixed}): For two simple graphs $X$ and $Y$ with the same (finite) number of vertices, their \textit{friends-and-strangers graph} $\FS(X,Y)$ is defined as:
\vspace*{-0.8em}
\begin{itemize}
    \setlength\itemsep{-0.4em}
    \item The vertex set $V(\FS(X,Y))$ of $\FS(X,Y)$ consists of all bijections $\sigma: V(X)\longrightarrow V(Y).$
    \item Two bijections $\sigma, \tau \in V(\FS(X,Y))$ are adjacent if and only if there exist two distinct vertices $u',v'\in V(X)$ such that:
    \vspace*{-0.8em}
    \begin{itemize}
    \setlength\itemsep{-0.4em}
        \item $u'$ and $v'$ are adjacent in $X,$
        \item $\sigma(u')$ and $\sigma(v')$ are adjacent in $Y,$
        \item $\sigma(u') = \tau(v')$ and 
        $\sigma(v') = \tau(u'),$
        \item $\sigma(w') = \tau(w')$ for all $w'\in V(X)\backslash \{u',v'\}.$
    \end{itemize}
\end{itemize}
\vspace*{-0.8em}
For two adjacent bijections $\sigma$ and $\tau$ that satisfy the above properties, we call the operation transforming $\sigma$ into $\tau$ (i.e., the transposition of $\sigma(u')$ and $\sigma(v')$) an \textit{$(X,Y)$-friendly swap}.\\

To illustrate the above definition, we interpret $X$ as a graph of $n$ positions, where an edge indicates that two positions are adjacent. $Y$ is interpreted as the friendship graph on a group of $n$ people. At each position $p',$ there is one person $\sigma(p').$ At any given point of time, two people are allowed to exchange positions (perform an $(X,Y)$-friendly swap) if and only if they are friends and are located in adjacent positions. We want to understand the behavior of this process. Some natural questions to consider are:
What is the set of configurations to which people can rearrange? In particular, can they rearrange to any given configuration? How long will it take them to do so? What is the long-term behavior of this process if at any given point of time, two randomly selected friends on adjacent positions perform a swap?\\

A more familiar example is the 15-puzzle, in which 15 blocks numbered 1 through 15 and an empty cell are placed on a $4\times 4$ grid. One is allowed to exchange the empty cell with any given block adjacent to it. In terms of friends-and-strangers graphs, the game is represented by a $4\times 4$ grid of positions $X\cong\mathsf{Grid}_{4\times 4}$ and a 16-vertex star graph $Y\cong \starreg_{16}$ with center the empty cell. Other ways to generalize the 15-puzzle besides friends-and-strangers graphs have also been studied \cite{rotatingpuzzle}.\\

Most known results on friends-and-strangers graphs address the first two questions raised above: "\textit{What is the set of configurations to which people can rearrange?}" and 
"\textit{In particular, can they rearrange to any given configuration?}"
Formally, these questions ask about the structure of the connected components of $\FS(X,Y)$ and, in particular, whether $\FS(X,Y)$ is connected. The results in literature regarding these questions roughly fall in three categories:
\vspace*{-0.4em}
\begin{itemize}
    \setlength\itemsep{-0.4em}
    \item \textbf{Concrete Structure:} If one or both of the graphs have a concrete structure, such as being a star graph \cite{Wilson} or a cycle or path \cite{FSFixed}, \cite{pathpath}.
    \item \textbf{Random Structure:} If both graphs $X$ and $Y$ are randomly generated Erdős–Rényi random graphs \cite{FSTypicalExtremal}.
    \item \textbf{Extremal Structure:} If the graphs $X$ and $Y$ have an extremal structure, such as a lower bound on their respective minimum degrees \cite{FSTypicalExtremal}. 
\end{itemize}

The current paper addresses questions in the third category. 

\subsection{Main Results}
In \cite{FSTypicalExtremal}, Alon, Defant, and Kravitz pose two related questions regrading the connected components of $\FS(X,Y)$ for $X$ and $Y$ with lower-bounded minimum degrees. The two cases they consider are, respectively, when $X$ and $Y$ are arbitrary graphs and when they are bipartite graphs. Studying bipartite graphs separately is natural, as will become apparent from Proposition 2.2.
\subsubsection{Arbitrary Graphs}

In \cite{FSTypicalExtremal}, the authors ask the following extremal question:\\

\noindent
\textbf{Problem 1.2 }(\cite{FSTypicalExtremal}):\textit{ Denote by $\delta(G)$ the minimum degree of a graph $G.$ Let $n$ be a natural number. What is the smallest natural number $d_n$ such that for any two graphs $X$ and $Y$ on $n$ vertices that satisfy $\delta(X)\ge d_n$ and $\delta(Y)\ge d_n,$ the graph $\FS(X,Y)$ is connected?}\\

In \cite[Theorem 1.3]{FSTypicalExtremal}, Alon, Defant, and Kravitz prove that 
$\frac{3n}{5}-2\le d_n$ holds for all $n,$ and $d_n\le \frac{9n}{14}+2$ holds when $n \ge 16.$ They conjecture \cite[Conjecture 7.3]{FSTypicalExtremal} that the true value of $d_n$ is $d_n = \frac{3n}{5}+O(1).$ As a corollary of one of our main results in the current paper --- Theorem 1.4 --- we find that $d_n \le \lceil 3n/5\rceil.$ Combined with the lower bound of $\frac{3n}{5}-2$ given in \cite[Theorem 1.3]{FSTypicalExtremal}, this settles the conjecture.\\

In the current paper, we address an asymmetric version of Problem 1.2 suggested in \cite{FSTypicalExtremal} which does not require the lower bounds on $\delta(X)$ and $\delta(Y)$ to be the same:\\

\noindent
\textbf{Problem 1.3} (\cite{FSTypicalExtremal}):\textit{ Find (sufficient and necessary) conditions on the pairs $(\delta_1(n),\delta_2(n))$ which guarantee that for any two connected $n$-vertex graphs $X$ and $Y$ such that $\delta(X)\ge \delta_1(n)$ and $\delta(Y)\ge \delta_2(n),$ the graph
$\FS(X,Y)$ is connected.}\\

Our main result regarding this problem is the following theorem:\\

\noindent
\textbf{Theorem 1.4}:\textit{
Suppose that $X$ and $Y$ are two graphs on $n\ge 6$ vertices satisfying:
\vspace*{-0.8em}
\begin{itemize}
    \setlength\itemsep{-0.4em}
    \item $\delta(X)>n/2, \delta(Y)>n/2,$
    \item $2\min(\delta(X),\delta(Y))+3\max(\delta(X), \delta(Y))\ge 3n.$
\end{itemize}
\vspace*{-0.8em}
Then $\FS(X,Y)$ is connected.}\\

An immediate corollary is that the value $d_n$ in Problem 1.2 satisfies $d_n \le \lceil 3n/5\rceil.$ As already mentioned, this resolves the conjecture \cite[Conjecture 7.3]{FSTypicalExtremal}.\\

We also prove another result of the same type as in Theorem 1.4, which gives weaker bounds, but also applies to the case when one of the graphs has degree at most $n/2:$\\

\noindent
\textbf{Theorem 1.5}:\textit{
Suppose that $X$ and $Y$ are two graphs on $n$ vertices satisfying:
\vspace*{-0.8em}
\begin{itemize}
    \setlength\itemsep{-0.4em}
    \item $X$ and $Y$ are both connected,
    \item $\min(\delta(X), \delta(Y))+2\max(\delta(X), \delta(Y))\ge 2n.$
\end{itemize}
\vspace*{-0.6em}
Then $\FS(X,Y)$ is connected.}\\

The proof of Theorem 1.5 is much simpler than the proof of Theorem 1.4 but still illustrates some of the essential ideas. For that reason, we prove Theorem 1.5 in Section 3, before we prove Theorem 1.4 in Section 4.\\

Our last result regarding Problem 1.3 focuses on lower bounds. It is a generalization of the lower bound construction in \cite[Theorem 1.3.]{FSTypicalExtremal} by Alon, Defant, and Kravitz:\\

\noindent
\textbf{Proposition 1.6}:\textit{ Suppose that $n\ge k\ge 5$ are integers. Then there exist connected graphs $X$ and $Y$ such that $\delta(X)\ge \frac{3n}{k}-4,\delta(Y)\ge \frac{(k-2)n}{k}-3$ and $\FS(X,Y)$ is disconnected.}\\

\begin{figure}[!htb]
\begin{center}
   \includegraphics[width = 
   .8\linewidth]{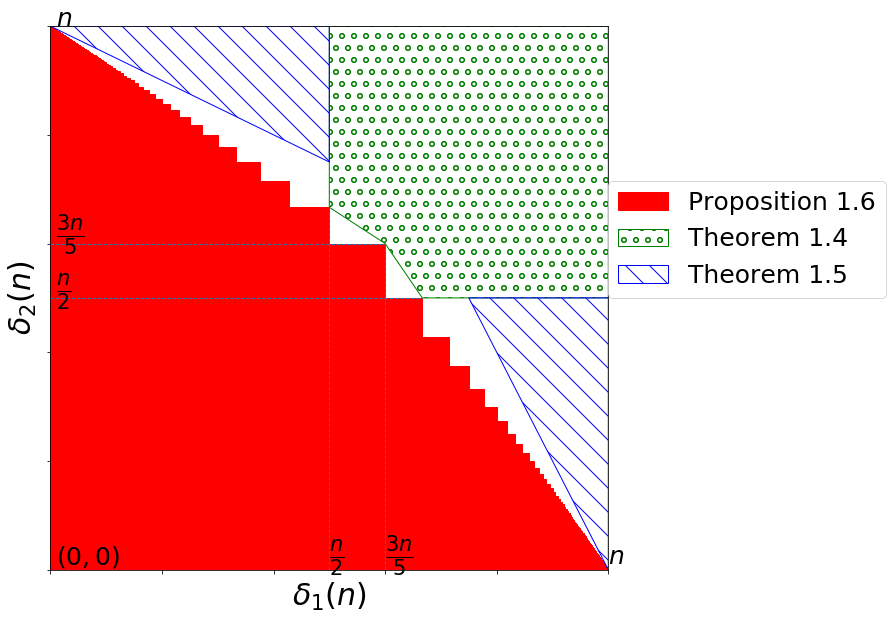}\caption{A visual representation of the results given by Theorem 1.4, Theorem 1.5, and Proposition 1.6 regarding Problem 1.3. In solid red are pairs $(\delta_1(n),\delta_2(n))$ which do not guarantee that $\FS(X,Y)$ is connected. In dotted green and striped blue are pairs that guarantee $\FS(X,Y)$ is connected. White portions are still unresolved.
   \textbf{Remark: }Additive constants are omitted in the diagram.}\label{fig:no_conj_diagram}
\end{center}
\end{figure}

The degrees $(\delta(X),\delta(Y))$ constructed in Proposition 1.6 lie strictly below the line in Theorem 1.5 parametrized by $\min(\delta(X), \delta(Y))+2\max(\delta(X), \delta(Y))= 2n$ and, thus, don't provide examples in which Theorem 1.5 is tight. However, they provide examples very close to the line $2\min(\delta(X),\delta(Y)+3\max(\delta(X), \delta(Y))= 3n$ from Theorem 1.4 --- as in \autoref{fig:no_conj_diagram}. This motivates Conjecture 8.1 extending Theorem 1.4 to the more general case when one of $X$ and $Y$ potentially has minimum degree at most $n/2$ but is still connected.

\newpage
\subsubsection{Bipartite Graphs}

In \cite{FSTypicalExtremal}, the authors also consider a version of Problem 1.2 when $X$ and $Y$ are taken to be edge-subgraphs of $K_{r,r},$ where $K_{r,r}$ denotes
the complete bipartite graph on two parts of size $r.$ As already mentioned, it will become clear from Proposition 2.2 why one asks for two connected components rather than one in this setting.\\ 

\noindent
\textbf{Problem 1.7 }(\cite{FSTypicalExtremal}):\textit{ What is the smallest natural number $d_{r,r},$ such that for any two edge-subgraphs $X$ and $Y$ of $K_{r,r}$ that satisfy $\delta(X)\ge d_{r,r}$ and $\delta(Y)\ge d_{r,r},$ the graph $\FS(X,Y)$ has exactly two connected components?}\\

They give a complete (up to an additive constant) answer to this question by proving the bounds $\displaystyle
\Big{\lceil} \frac{3r+1}{4}\Big{\rceil} \le d_{r,r}\le 
\Big{\lceil} \frac{3r+2}{4}\Big{\rceil}
$ in \cite[Theorem 1.4]{FSTypicalExtremal}.\\

Alon, Defant, and Kravitz conjecture \cite[Conjecture 7.4]{FSTypicalExtremal} that the true value of $d_{r,r}$ is $\lceil \frac{3r+1}{4}\rceil.$ They also pose another problem related to Problem 1.7:\\

\noindent
\textbf{Problem 1.8 }(\cite{FSTypicalExtremal}):\textit{ What is the smallest natural number $d^*_{r,r},$ such that for any edge-subgraph $X$ of $K_{r,r}$ that satisfies $\delta(X)\ge d^*_{r,r},$ the graph $\FS(X,K_{r,r})$ has exactly two connected components?}\\

In analogy to Problem 1.3, we ask about an asymmetric version of Problem 1.7. The asymmetric version, in particular, also addresses Problem 1.8.\\

\noindent
\textbf{Problem 1.9}:\textit{ Find (sufficient and necessary) conditions on the pairs $(\delta_1(r,r),\delta_2(r,r))$ which guarantee that for any two edge-subgraphs $X$ and $Y$ of $K_{r,r}$  such that $\delta(X)\ge \delta_1(r,r)$ and $\delta(Y)\ge \delta_2(r,r),$ the graph
$\FS(X,Y)$ has exactly two connected components.}\\

Our main result regarding Problem 1.9 is the following theorem:\\

\noindent
\textbf{Theorem 1.10}: \textit{Let $r\ge 2,$ and let $X$ and $Y$ be edge-subgraphs of $K_{r,r}$ such that $$\delta(X)+\delta(Y)\ge 3r/2+1.$$ Then $\FS(X,Y)$ has exactly two connected components}.\\

We obtain a complementary lower-bound result by extending the construction in \cite[Theorem 1.4]{FSTypicalExtremal}:\\

\noindent
\textbf{Theorem 1.11}: \textit{ Let $r\ge 2$ and $\delta_1, \delta_2$ be non-negative integers satisfying
$$\delta_1 + \delta_2 = \Big{\lfloor} \frac{3r}{2}\Big{\rfloor} ,\; \delta_1 \le r,\;  \delta_2\le r.$$ 
Then there exist two edge-subgraphs $X$ and $Y$ of $K_{r,r}$ such that $\delta(X) = \delta_1, \delta(Y) = \delta_2,$ and $\FS(X,Y)$ has more than two connected components.}\\

Theorem 1.11 certifies that the condition $\delta(X)+\delta(Y)\ge 3r/2+1$ is tight when $r$ is even and is off by at most 1 when $r$ is odd.\\

Setting $\delta_2 = r$ in Theorems 1.10 and 1.11,  we obtain the following (almost complete) answer to Problem 1.8:\\

\noindent
\textbf{Corollary 1.12}:\textit{ The number $d_{r,r}^*$ from Problem 1.8 satisfies 
\begin{equation*}
  d^*_{r,r} =
    \begin{cases}
      r/2+1 & \text{ when $r$ is even.}\\
      \lceil r/2\rceil \text{ or } \lceil r/2\rceil+1 & \text{ when $r$ is odd.}
    \end{cases}       
\end{equation*}
}

Finally, we note that our proof of Theorem 1.10 also provides an alternative argument for the upper bound $\displaystyle d_{r,r}\le 
\Big{\lceil} \frac{3r+2}{4}\Big{\rceil}
$ in \cite[Theorem 1.4]{FSTypicalExtremal} since $\displaystyle 2
\Big{\lceil} \frac{3r+2}{4}\Big{\rceil}\ge 3r/2 + 1.
$

\section{Preliminaries}

\subsection{Terminology and Notation}
Throughout, we will use the following terminology and notation. For a simple graph $X:$
\vspace*{-0.6em}
\begin{itemize}
    \setlength\itemsep{-0.4em}
    \item Denote by $V(X)$ the set of vertices and by $E(X)$ the set of edges of $X.$
    \item For a vertex $x$ of $X,$ denote by $N(x)$ the \textit{open neighborhood} of $x$ in $X,$  which is given by the set of neighbors of $x.$ Denote by $N[x]$ the \textit{closed neighborhood} of $x$ in $X,$ given by $N[x] = N(x)\cup\{x\}.$
    \item For any $A\subseteq V(X),$ denote by $X|_A$ the induced subgraph of $X$ on the vertex set $A.$
    \item Denote by $\delta(X)$ the minimum degree of $X.$
\end{itemize}

\begin{figure}[!htb]
\begin{center}
   \includegraphics[width = 0.4\linewidth]{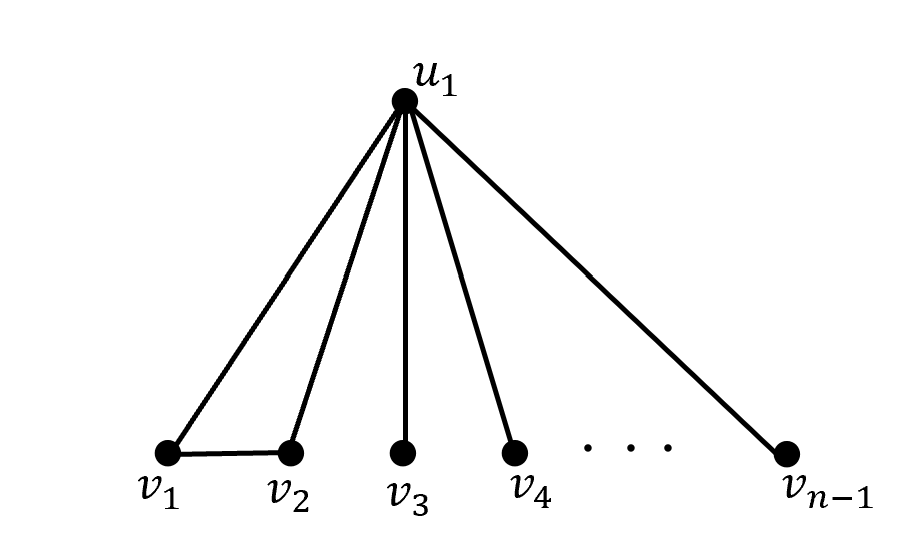}\caption{$\starplus_n$ Graph.}\label{fig:starplus}
\end{center}
\end{figure}

We also introduce notation for the following special graphs: 
\vspace*{-0.5em}
\begin{itemize}
    \setlength\itemsep{-0.4em}
    \item $K_n$ is the complete graph on $n$ vertices.
    \item $K_{n,m}$ is the complete bipartite graph with parts of sizes $n$ and $m.$
    \item $\starreg_n$ is $K_{1,n-1}.$
    \item $\starplus_n$ is a star with $n$ vertices and an extra edge, as in \autoref{fig:starplus}.
\end{itemize}
\vspace*{-0.4em}

Finally, we introduce notation specific to friends-and-strangers graphs. Denote an $(X,Y)$-friendly swap that transforms a bijection $\sigma$ into the bijection $\tau = (u,v)\circ\sigma$ by $uv$ (recall that adjacent bijections in $\FS(X,Y)$ differ by a transposition in $V(Y)$). Similarly, we can have a sequence of more than one $(X,Y)$-friendly swap. For example, consider the graphs $X$ (in \autoref{fig:seqX}) and $Y$ (in \autoref{fig:seqY}) below. The bijection $\sigma: V(X)\longrightarrow V(Y),$ where $\sigma(1) = a, \sigma(2) = b, \sigma(3) = c, \sigma(4) = d,$ is also given.

\begin{figure}[!htb]
\minipage{0.45\textwidth}
\begin{center}
  \includegraphics[width=0.5\linewidth]{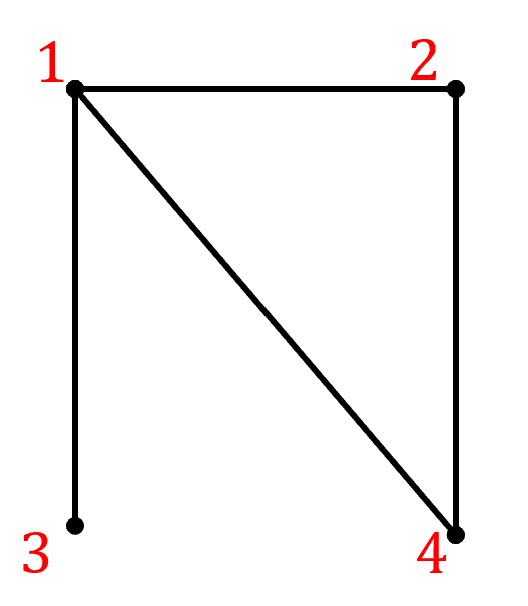}
  \caption{$X$ Graph}\label{fig:seqX}
  \end{center}
\endminipage\hfill
\minipage{0.45\textwidth}%
\begin{center}
  \includegraphics[width=0.5\linewidth]{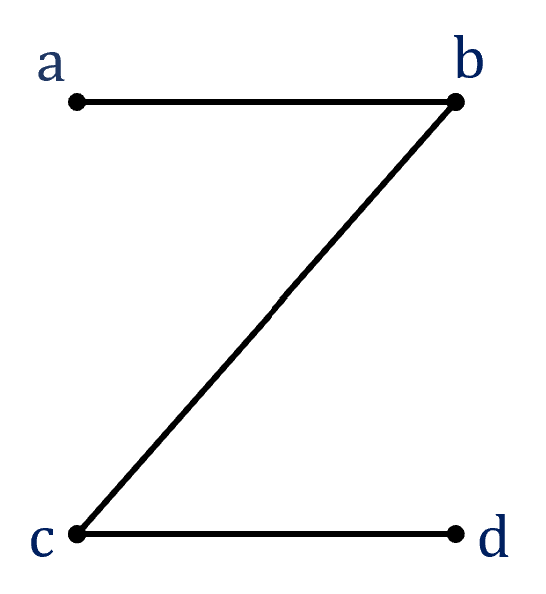}
  \end{center}
  \caption{$Y$ Graph}\label{fig:seqY}
\endminipage
\end{figure}

\begin{figure}[!htb]
\begin{center}
   \includegraphics[width = 0.8\linewidth]{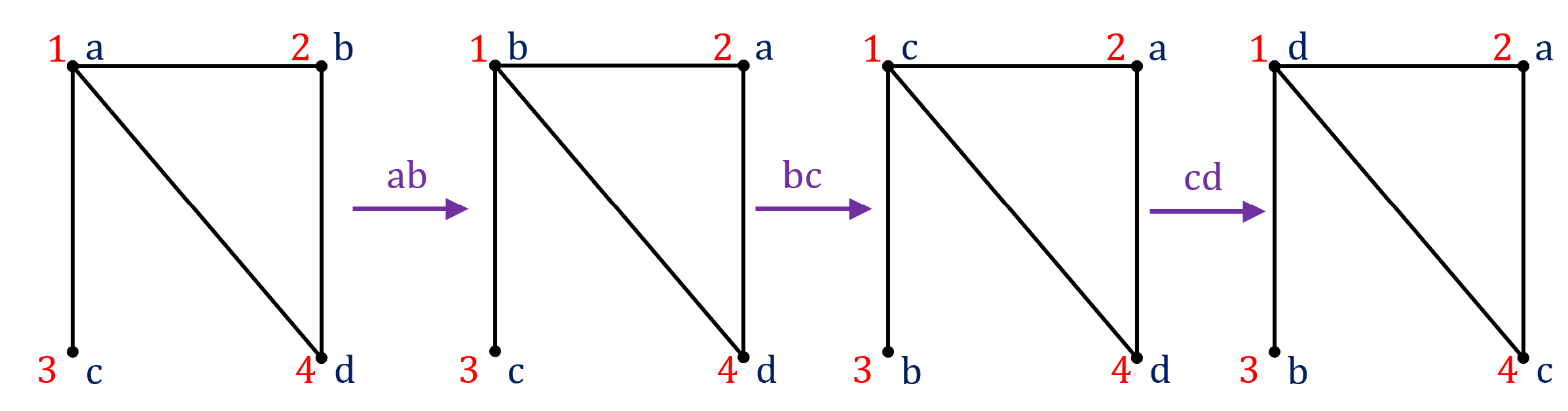}\caption{The sequence $\Sigma = ab,bc,cd$ applied to the graphs in \autoref{fig:seqX} and \autoref{fig:seqY}.}\label{fig:seqProcess}
\end{center}
\end{figure}

Then the sequence $\Sigma = ab, bc, cd$ transforms $\sigma$ into the bijection $\tau$ given by $\tau(1) = d,\tau(2) = a,\tau(3) = b,$ and $\tau(4) = d$ as in \autoref{fig:seqProcess}.

Observe that if the sequence of $(X,Y)$-friendly swaps $\Sigma$ transforms the bijection $\sigma$ into the bijection $\tau,$ then the reverse sequence, which we denote by $\rev(\Sigma)$ throughout the rest of the paper, will also be a sequence of $(X,Y)$-friendly swaps and will transform $\tau$ back into $\sigma.$\\

Finally, note that if $X'$ is a subgraph of $X$ and $Y'$ is a subgraph of $Y,$ then any sequence $\Sigma'$ of $(X',Y')$-friendly swaps is also a sequence of $(X,Y)$-friendly swaps. Furthermore, if  $x$ is a vertex in $V(X)\backslash V(X')$ and the sequence $\Sigma'$ of $(X',Y')$-friendly swaps transforms the bijection\\ $\sigma: V(X)\longrightarrow V(Y)$ into the bijection $\tau: V(X)\longrightarrow V(Y),$ then $\sigma(x) = \tau(x).$ 

\subsection{Related Results on Friends-And-Strangers Graphs}
Here we summarize previous work on friends-and-strangers graphs relevant to the current paper. 

\subsubsection{General Properties}
We begin with a proposition showing that $X$ and $Y$ play symmetric roles in $\FS(X,Y):$\\

\noindent
\textbf{Proposition 2.1 }(\cite{FSFixed}): \textit{If $X$ and $Y$ are $n$-vertex graphs, then $\FS(X,Y)$ and $\FS(Y,X)$ are isomorphic.}\\

We also discuss two obstructions to the connectivity of $\FS(X,Y).$
A simple one is if either $X$ or $Y$ is disconnected. For that reason, in Theorem 1.5 we explicitly require that $X$ and $Y$ are connected (the condition is implied in Theorem 1.4 by $\delta(X)>n/2,\delta(Y)>n/2$).
Another, less trivial, obstruction to the connectivity of $\FS(X,Y)$ is given by the following proposition in \cite{FSFixed}. An equivalent to this statement for the special case of the 15-puzzle has been known at least since 1879 \cite{15puz}.\\

\noindent
\textbf{Proposition 2.2 }(\cite[Proposition 2.7]{FSFixed}): \textit{If $X$ and $Y$ are bipartite graphs, each on $n \ge 3$ vertices, then $\FS(X,Y)$ is disconnected.
In particular, if the partite sets of $X$ are $A_X\sqcup B_X$  and the partite sets of $Y$ are $A_Y\sqcup B_Y,$ no two bijections $\sigma:V(X)\longrightarrow V(Y)$ and $\tau:V(X)\longrightarrow V(Y)$ for which $$sign(\tau^{-1}\circ\sigma)\not \equiv
|\sigma(A_X)\cap A_Y| - |\tau(A_X)\cap A_Y| \pmod{2}
$$
are in the same connected component of $\FS(X,Y).$ As usual, $sign$ is the non-trivial homomorphism from $S_n$ to the cyclic group $\{0,1\}.$
}\\

Proposition 2.2 sets a natural barrier to the connectivity of $\FS(X,Y)$ when $X$ and $Y$ are taken from a bipartite family of graphs such as (random) subgraphs of $K_{r,r}.$ Thus, in such settings it is natural to ask for conditions guaranteeing that $\FS(X,Y)$ has exactly two connected components, which is the fewest possible. We study such existence of exactly 2 connected components in Theorems 1.10 and 1.11, which generalize \cite[Theorem 1.4]{FSTypicalExtremal}. In proving the existence of exactly two connected components, we will also use the following special case of \cite[Proposition 2.6]{FSTypicalExtremal}.\\

\noindent
\textbf{Proposition 2.3 }(\cite{FSTypicalExtremal}): \textit{Let $r\ge 2.$ Then $\FS(K_{r,r}, K_{r,r})$ has exactly two connected components.}

\subsubsection{Star Graphs}
The special case when $X$ is a star graph, which generalizes the 15-puzzle, was studied by Wilson in 1974 \cite{Wilson}. This case is of significant importance to the extremal setting we are studying as large star subgraphs appear naturally when one has a lower bound on the minimum degree. Indeed, in a graph $X,$ every vertex $x$ is the center of a star subgraph with at least $\delta(X)+1$ vertices, namely the vertices in the closed neighborhood $N[x].$ Identifying appropriate (parts of) closed neighborhoods to which one can apply Wilson's Theorem is at the heart of the proofs of many results on friends-and-strangers graphs, such as the ones given in \cite[Theorem 1.3]{FSTypicalExtremal}, Theorem 1.4, and Theorem 1.5. To introduce Wilson's result on star graphs, we first make the following definitions:\\

\noindent
\textbf{Definition 2.4}: Let $X$ be a connected graph. A \textit{cut-vertex} is a vertex $v\in V(X)$ such that $X|_{V(X)\backslash \{v\}}$ is not connected.\\

\noindent
\textbf{Definition 2.5}: A graph $X$ that is connected and has no cut-vertex is \textit{biconnected}.\\

\begin{figure}[!htb]
\begin{center}
   \includegraphics[width = 0.22\linewidth]{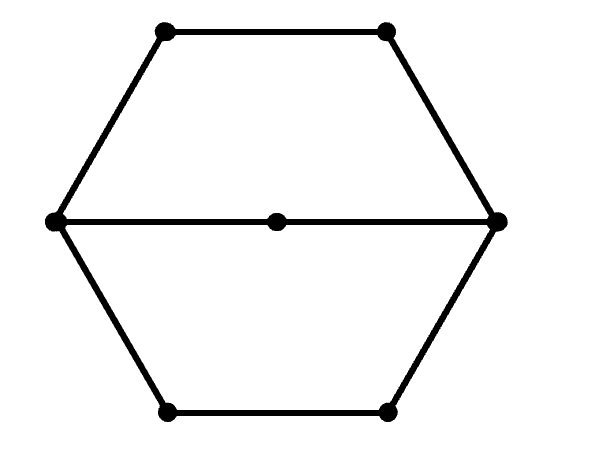}\caption{The Graph $\theta_0.$}\label{fig:theta0}
\end{center}
\end{figure}

\noindent
\textbf{Theorem 2.6 }(\cite[Theorem 1]{Wilson}): \textit{Suppose that $Y$ is a graph on $n$ vertices satisfying the four properties:
\vspace*{-.8cm}
\begin{itemize}
    \setlength\itemsep{-0.5em}
    \item $Y$ is biconnected, 
    \item $Y$ is not bipartite,
    \item $Y$ is not isomorphic to a cycle graph,
    \item $Y$ is not isomorphic to the graph on $7$ vertices denoted by $\theta_0$ as in \autoref{fig:theta0}.
\end{itemize}
\vspace*{-0.5em}
Then $\FS(\starreg_n, Y)$ is connected.
}

From now on, we we will call graphs satisfying the 4 properties listed in Theorem 2.6 \textit{Wilsonian}. The following lemma due to Alon, Defant, and Kravitz provides a useful tool for identifying Wilsonian graphs and is another main technique in the upper-bound proofs of \cite[Theorem 1.3.]{FSTypicalExtremal}, and Theorems 1.4 and 1.5.\\

\noindent
\textbf{Lemma 2.7} (\cite[Lemma 2.4]{FSTypicalExtremal}): \textit{Suppose that $X$ is a graph on $n$ vertices such that $\delta(X)>n/2.$ Then $X$ is Wilsonian.} \\

Of particular importance to the proof of Theorem 1.4 is also the following corollary of Theorem 2.6 and Proposition 2.2 due to Defant and Kravitz.\\

\noindent
\textbf{Theorem 2.8 }(\cite[Remark 2.8]{FSFixed}): \textit{Suppose that $Y$ is a graph on $n\ge 3$ vertices satisfying the three properties:
\vspace*{-0.4em}
\begin{itemize}
    \setlength\itemsep{-0.4em}
    \item $Y$ is biconnected, 
    \item $Y$ is not isomorphic to a cycle graph on at least 4 vertices,
    \item $Y$ is not isomorphic to $\theta_0$ (\autoref{fig:theta0}).
\end{itemize}
\vspace*{-0.4em}
Then $\FS(\starplus_n, Y)$ is connected.}\\

As graphs satisfying the three conditions listed in Theorem 2.8 will appear frequently in the upcoming proofs, we will call them \textit{almost-Wilsonian} for brevity.

\subsubsection{Exchangeable Pairs}
The last concept we need to introduce before beginning the proofs of Theorems 1.4 and 1.5 is that of exchangeable pairs of vertices, again due to Alon, Defant, and Kravitz.\\

\noindent
\textbf{Definition 2.9 }(\cite{FSTypicalExtremal}): Suppose that $X$ and $Y$ are two $n$-vertex graphs. A bijection $\sigma: V(X)\longrightarrow V(Y)$ and two distinct vertices $u,v\in V(Y)$ are given. We say that \textit{$u$ and $v$ are  $(X,Y)$-exchangeable from $\sigma$} if there exists a sequence of $(X,Y)$-friendly swaps that transforms $\sigma$ into $(u,v)\circ \sigma.$\\
\textbf{Remark:} Equivalently, we can define $u$ and $v$ being $(X,Y)$-exchangeable from $\sigma$ if $\sigma$ and $(u,v)\circ \sigma$ are in the same connected component of $\FS(X,Y).$\\

The following observations in \cite{FSTypicalExtremal} illustrate how the concept of exchangeable pairs is useful in studying the connectivity of a graph $\FS(X,Y):$\\

\noindent
\textbf{Proposition 2.10 }(\cite[Proposition 2.8]{FSTypicalExtremal}): \textit{Let $X,Y$ and $\widetilde{Y}$
be $n$-vertex graphs such that $Y$ is an edge-subgraph of $\widetilde{Y}.$ Suppose that for every edge $\{u,v\}$ of $\widetilde{Y}$ and every bijection $\sigma$ satisfying $(\sigma^{-1}(u), \sigma^{-1}(v))\in E(X),$ the vertices $u$ and $v$ are $(X,Y)$-exchangeable from $\sigma.$ Then the connected components of $\FS(X,Y)$ and the connected components of $\FS(X, \widetilde{Y})$ have the same vertex sets. In particular, the number of connected components of $\FS(X, \widetilde{Y})$ is equal to the number of connected components of $\FS(X,Y).$}\\

Of particular importance to the current paper is the following corollary of Proposition 2.10 in \cite{FSTypicalExtremal}:\\

\noindent
\textbf{Lemma 2.11} (\cite[Lemma 2.9]{FSTypicalExtremal}): \textit{Let $X$ and $Y$ be $n$-vertex graphs with $X$ connected. Suppose that for all distinct vertices $u,v\in V(Y)$ and bijections $\sigma: V(X)\longrightarrow V(Y)$ such that $(\sigma^{-1}(u), \sigma^{-1}(v))\in E(X),$ the vertices $u$ and $v$ are $(X,Y)$-exchangeable from $\sigma.$ Then $\FS(X,Y)$ is connected.}\\

The proofs of \cite[Theorem 1.3]{FSTypicalExtremal}, Theorem 1.4, and Theorem 1.5 use Lemma 2.11 as follows. One first chooses an arbitrary bijection $\sigma: V(X)\longrightarrow V(Y)$ and two vertices $u,v\in V(Y)$ such that  $(\sigma^{-1}(u), \sigma^{-1}(v))\in E(X).$ The goal is to prove that
$u$ and $v$ are $(X,Y)$-exchangeable from $\sigma.$\\

\begin{figure}[!htb]
\begin{center}
   \includegraphics[width = .9\linewidth]{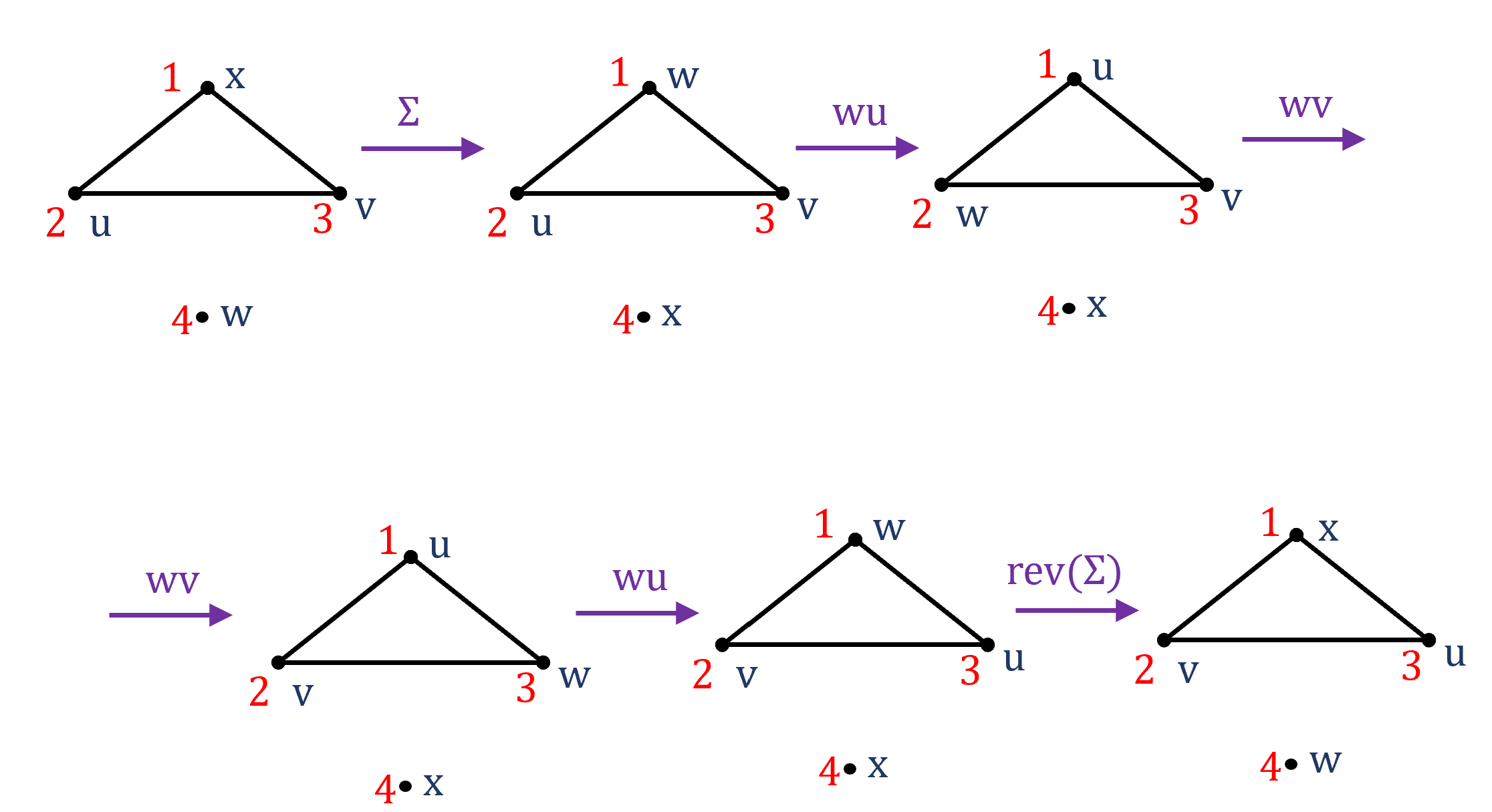}\caption{The sequence $\Sigma,wu,wv,wu,\rev(\Sigma)$ applied to the permutation $\sigma,$ where $\sigma(1) =x,$ $\sigma(2) = u, \sigma(3) = v,$ and $\sigma(4) = w.$ As a result, $u$ and $v$ are exchanged.}\label{fig:viaNeighbors}
\end{center}
\end{figure}

In the proofs of Theorems 1.4 and 1.5 in the current paper, several techniques for proving that the so chosen vertices $u$ and $v$ are $(X,Y)$-exchangeable from $\sigma$ reappear.
\begin{enumerate}
    \setlength\itemsep{-0.4em}
    \item \underline{Via trivial swaps:} Here are two examples. A first example is if $u$ and $v$ are adjacent. In that case, $uv$ is an $(X,Y)$-friendly swap that exchanges them. A second example is if there exists a vertex $w\in V(Y)$ such that $(u,w)\in E(Y), (v,w)\in E(Y),$ $(\sigma^{-1}(u), \sigma^{-1}(w))\in E(X),$ and $(\sigma^{-1}(v), \sigma^{-1}(w))\in E(X).$ In that case, the sequence of $(X,Y)$-friendly swaps $wu,wv,wu$ exchanges $u$ and $v$ as in \autoref{fig:viaTrivial}.
    \begin{figure}[!htb]
    \begin{center}
   \includegraphics[width = \linewidth]{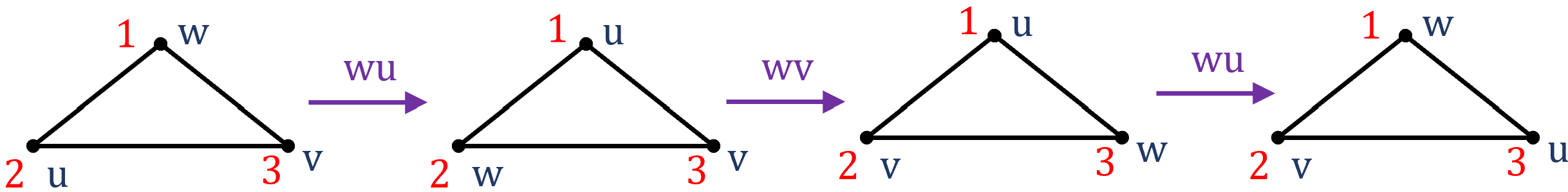}\caption{The sequence $wu, wv,wu$ of $(X,Y)$-friendly swaps applied to the bijection $\sigma$ from the second example in 1., which is given by $\sigma^{-1}(w) =1,
   \sigma^{-1}(u) =2, \sigma^{-1}(v) =3.$}\label{fig:viaTrivial}
\end{center}
\end{figure}
    \item \underline{Via Theorem 2.6:} One identifies a subgraph $X_1$ of $X$ and a subgraph $Y_1$ of $Y$ which satisfy the following two properties:
    \vspace*{-0.6em}
    \begin{itemize}
    \setlength\itemsep{-0.4em}
        \item $\sigma^{-1}(u), \sigma^{-1}(v)\in V(X_1)$ and $u,v\in V(Y_1),$ and
        \item Either $X_1$ contains a spanning star and $Y_1$ is Wilsonian, or $Y_1$ contains a spanning star and $X_1$ is Wilsonian. 
    \end{itemize} 
    \vspace*{-0.6em}
    Then one invokes Theorem 2.6 for $X_1$ and $Y_1.$
    \vspace*{0.4em}
    \item \underline{Via Theorem 2.8:} One identifies a subgraph $X_1$ of $X$ and a subgraph $Y_1$ of $Y,$ which satisfy the following two properties:
    \begin{itemize}
    \setlength\itemsep{-0.4em}
        \item $\sigma^{-1}(u), \sigma^{-1}(v)\in V(X_1)$ and $u,v\in V(Y_1),$ and
        \item Either $X_1$ contains a spanning $\starplus$ graph and $Y_1$ is almost-Wilsonian, or $Y_1$ contains a spanning $\starplus$ graph and $X_1$ is almost-Wilsonian.  
    \end{itemize}
    \vspace*{-0.6em}
    Then one invokes Theorem 2.8 for the graphs $X_1$ and $Y_1.$
    \vspace*{0.4em}
    \item \underline{Via common neighbors:} One identifies two vertices $x$ and $w$ of $Y$ and a sequence $\Sigma,$ which satisfy the following three properties:
    \vspace*{-0.4em}
    \begin{itemize}
        \item $(\sigma^{-1}(x),\sigma^{-1}(u))\in E(X), (\sigma^{-1}(x),\sigma^{-1}(v))\in E(X),$ and
        \item $(w,u)\in E(Y), (w,v)\in E(Y),$ and
        \item $\Sigma$ is a sequence of $(X,Y)$-friendly swaps that does not involve $u$ and $v$ and transforms $\sigma$ into a permutation $\tau$ satisfying $\tau^{-1}(w) = \sigma^{-1}(x).$
    \end{itemize}
    \vspace*{-0.4em}
     Having identified such $w,x$ and $\Sigma,$ one can readily check that the sequence of $(X,Y)$-friendly swaps $\Sigma, wu, wv, wu, \rev(\Sigma)$ exchanges $u$ and $v$ as shown in \autoref{fig:viaNeighbors}. 
\end{enumerate}

One major difference between the proof of the upper bound in \cite[Theorem 1.3]{FSTypicalExtremal} and the proof of the stronger upper bound in Theorem 1.4 in this current paper is the third idea of identifying almost-Wilsonian graphs in addition to Wilsonian graphs.

\section{Upper Bounds for Arbitrary Graphs: Theorem 1.5}
\noindent
\textbf{Theorem 1.5}:\textit{
Suppose that $X$ and $Y$ are two graphs on $n$ vertices satisfying:
\vspace*{-0.6em}
\begin{itemize}
    \setlength\itemsep{-0.4em}
    \item $X$ and $Y$ are both connected,
    \item $\min(\delta(X), \delta(Y))+2\max(\delta(X), \delta(Y))\ge 2n.$
\end{itemize}
\vspace*{-0.6em}
Then $\FS(X,Y)$ is connected.}\\

\noindent
\begin{proof}
Suppose that $X$ and $Y$ are two graphs satisfying the above properties. We can assume without loss of generality that $\delta(X)\le \delta(Y)$ as $\FS(X,Y)$ and $\FS(Y,X)$ are isomorphic by Proposition 2.1. Thus, the second condition simplifies to $\delta(X)+2\delta(Y)\ge 2n.$\\

Fix a bijection
$\sigma: V(X)\longrightarrow V(Y)$ and $u,v\in V(Y)$ such that 
$(\sigma^{-1}(u), \sigma^{-1}(v))\in E(X).$ We will prove that $u$ and $v$ are
$(X,Y)$-exchangeable from $\sigma.$ As $X$ is connected and the vertices $u,v$ were chosen arbitrarily, 
it will then follow from Lemma 2.11 that $\FS(X, Y)$ is connected as desired.\\

Consider the set $Q = \sigma(N[\sigma^{-1}(u)]).$ We will prove that $\delta(Y|_{Q})> \frac{1}{2}|Q|.$ Take any vertex $y \in Q.$ Its degree in $Y|_Q$ is
$$
|N(y)\cap Q| = |N(y)|+|Q| - |N(y)\cup Q|\ge |N(y)|+|Q|-n,
$$
so $\delta(Y|_Q)\ge |N(y)|+|Q|-n.$
It thus suffices to show that $|N(y)|+|Q|-n>\frac{1}{2}|Q|.$ This is equivalent to
$$
2|N(y)|+ |Q|>2n.
$$

However, $|N(y)|\ge \delta(Y)$ by the definition of $\delta,$ and $|Q| = |\sigma(N[\sigma^{-1}(u)])| =|N[\sigma^{-1}(u)]|\ge \delta(X)+1.$ Thus,
$$
2|N(y)|+ |Q|\ge 2\delta(Y)+\delta(X)+1\ge 2n + 1 >2n,
$$
as desired. We illustrate the sets and vertices considered in \autoref{fig:Thm1.6}.\\

\begin{figure}[!htb]
\begin{center}
   \includegraphics[width = \linewidth]{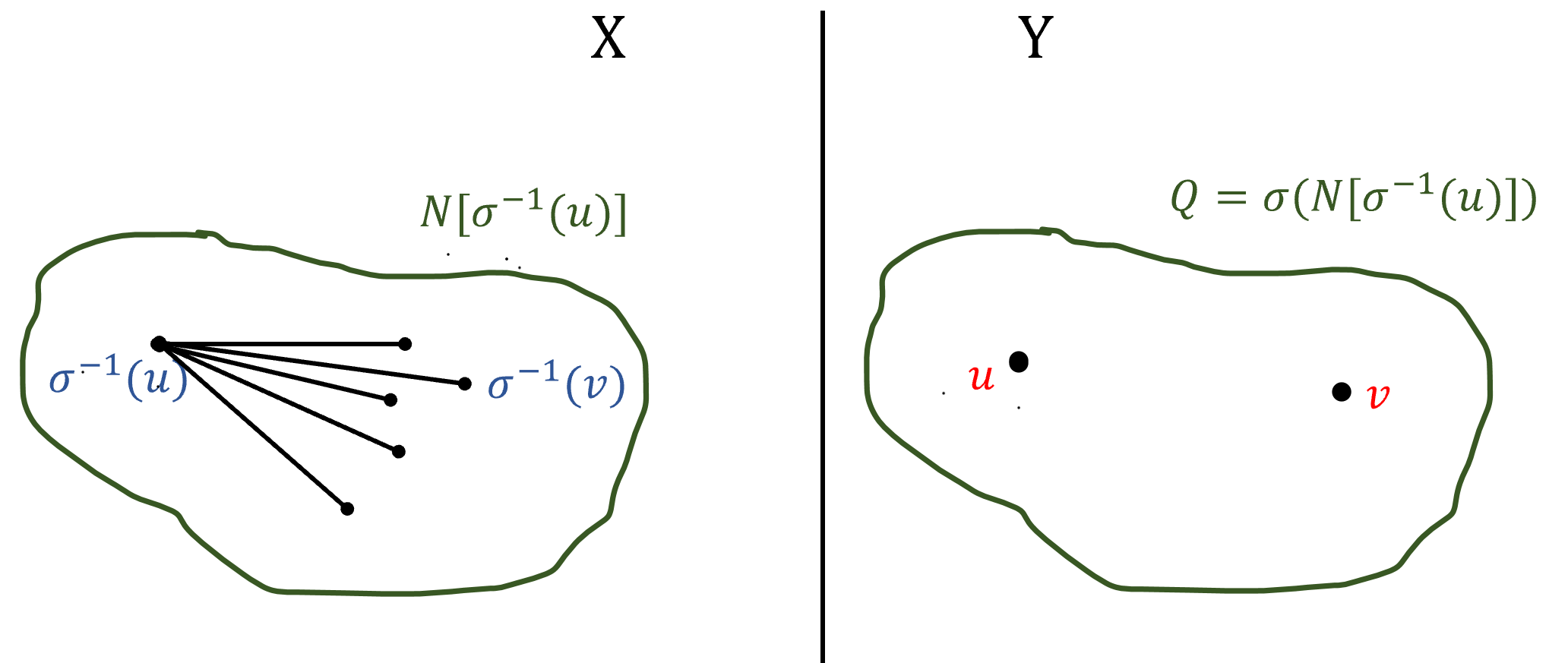}\caption{The sets $N[\sigma^{-1}(u)]$ and $Q = \sigma(N[\sigma^{-1}(u)]).$}\label{fig:Thm1.6}
\end{center}
\end{figure}

Since $\delta(Y|_{Q})>|Q|/2,$ it follows from Lemma 2.7 that $Y|_{Q}$ is Wilsonian. At the same time, 
$X|_{\sigma^{-1}(Q)} = X|_{N[\sigma^{-1}(u)]}$ clearly contains a spanning star with center $\sigma^{-1}(u).$ Thus, by Theorem 2.6, 
the graph $\FS(X|_{\sigma^{-1}(Q)}, Y|_Q)$ is connected. In particular, as $\sigma^{-1}(v)\in N(\sigma^{-1}(u))\subseteq \sigma^{-1}(Q),$ there exists a sequence $\Sigma$ of $(X|_{\sigma^{-1}(Q)}, Y|_Q)$-friendly swaps transforming the bijection $\sigma|_{\sigma^{-1}(Q)}$ into the bijection $(u,v)\circ \sigma|_{\sigma^{-1}(Q)}.$ When $\Sigma$ is viewed as sequence of $(X,Y)$-friendly swaps that do not involve vertices out of $Q,$ $\sigma$ is clearly transformed by $\Sigma$ into $(u,v)\circ\sigma.$ It follows that $u$ and $v$ are
$(X,Y)$-exchangeable from $\sigma,$ which completes the proof.
\end{proof}

\section{Upper Bounds for Arbitrary Graphs: Theorem 1.4}

\noindent
\textbf{Theorem 1.4}:\textit{
Suppose that $X$ and $Y$ are two graphs on $n\ge 6$ vertices satisfying:
\vspace*{-0.8em}
\begin{itemize}
    \setlength\itemsep{-0.4em}
    \item $\delta(X)>n/2, \delta(Y)>n/2,$
    \item $2\min(\delta(X), \delta(Y))+3\max(\delta(X),\delta(Y))\ge 3n.$
\end{itemize}
\vspace*{-0.6em}
Then $\FS(X,Y)$ is connected.}\\

We will split the proof into multiple sections.

\subsection{On the Number of Connected Components of Certain Subgraphs}
An important preliminary claim is the following lemma:\\

\noindent
\textbf{Lemma 4.1}: 
\textit{Let $G$ be a graph on $m$ vertices with minimum degree $\delta(G).$ Let $Q$ be any subset of $V(G)$ such that $|Q|\ge 5$ and $2|Q|+3\delta(G)\ge 3m +2$. Then $G|_Q$ has at most two connected components. Furthermore:\\
1) If $G|_Q$ has exactly two connected components, $F_1$ and $F_2,$ then both $F_1$ and $F_2$ are Wilsonian, and the following inequalities are satisfied: $$\delta(G)+1+|Q|-m\le |V(F_1)|\le m - \delta(G)-1,\; \; \delta(F_1)\ge \delta(G)+|Q|-m,$$ 
$$\delta(G)+1+|Q|-m\le |V(F_2)|\le m - \delta(G)-1,\; \; \delta(F_2)\ge\delta(G)+|Q|-m.$$
2) If $G|_Q$ has a single connected component $F,$ then one of the following holds:
\vspace*{-0.8em}
\begin{itemize}
\setlength\itemsep{-0.4em}
    \item $F$ is almost-Wilsonian.
    \item There exists a cut vertex $v$ such that $F|_{V(F)\backslash \{v\}}$ has exactly two connected components $F_1'$ and $F_2'.$ Furthermore, both $F_1', F_2'$ are Wilsonian, and the following inequalities are satisfied:
    $$\delta(G)+|Q|-m\le |V(F'_1)|\le m - \delta(G)-1,\; \; \delta(F'_1)\ge \delta(G)+|Q|-m-1,$$ 
$$\delta(G)+|Q|-m\le |V(F'_2)|\le m - \delta(G)-1,\; \; \delta(F'_2)\ge\delta(G)+|Q|-m-1.$$
Finally, if $|N(v)\cap V(F_1')|\ge 2,$ then the graph $F|_{V(F_1')\cup \{v\}}$ is Wilsonian, and likewise for $F_2'.$
\end{itemize}
}

\noindent
\begin{proof} We begin with a lower bound on $\delta(G|_Q).$ Take any vertex $q\in Q.$ Then its degree in $G|_{Q}$ is
$$
|N(q)\cap Q|\ge |N(q)| + |Q|-m\ge \delta(G)+|Q|-m, 
$$
so $\delta(G|_Q)\ge \delta(G)+|Q|-m.$
In particular, this means that the connected component of $q$ in $G|_Q$ has size at least $\delta(G)+|Q|-m+1.$\\
First, we will use this observation to prove that the number of connected components is at most 2.
Suppose, for the sake of contradiction, that $G|_Q$ has at least three connected components, $K_1, K_2,$ and $K_3.$ As above, $|V(K_i)|\ge \delta(G)+|Q|-m+1$ for each $i.$ Since $V(K_1), V(K_2),$ and  $V(K_3)$ are disjoint subsets of $Q:$
$$
|Q|\ge |V(K_1)|+|V(K_2)|+|V(K_3)|\ge 3(|Q|+\delta(G)+1-m)\Longrightarrow
$$
$$
3m-3 \ge 2|Q|+3\delta(G),
$$
which contradicts the assumed inequality  $2|Q|+3\delta(G)\ge 3m +2.$ To prove the rest of the statement, we separately consider the two cases based on the number of connected components of $G|_{Q}.$\\

\noindent
\textbf{Case I:} Suppose that $G|_{Q}$ has two connected components $F_1$ and $F_2.$ As $F_1$ is a connected component, we already know that
$|V(F_1)|\ge \delta(G)+|Q|-m+1$. Now, let $u$ be an arbitrary vertex in $F_2.$ Clearly $V(F_1)\cap N[u] = \emptyset.$ Thus,
$$
|V(F_1)|\le m - |N[u]| \le m - \delta(G)-1.
$$
Next, as $F_1$ is a connected component in $G|_Q,$ clearly $\delta(F_1)\ge \delta(G|_{Q}).$ It follows that
$$
\delta(F_1)\ge \delta(G|_Q)\ge \delta(G)+|Q|-m > (m - \delta(G)-1)/2\ge |V(F_1)|/2,
$$
again using the assumed inequality $2|Q|+3\delta(G)\ge 3m +2.$ Therefore, $F_1$ is Wilsonian by Lemma 2.6. So, $F_1$ satisfies all the conclusions of the lemma. The analogous reasoning shows that $F_2$ also does so. \\

\noindent
\textbf{Case II:} Suppose that $G|_Q$ has a single connected component $F.$ We consider two cases:\\
II.1) Suppose that $F$ is biconnected.
To show that $F$ is almost-Wilsonian, we simply need to show that it is not isomorphic to a cycle graph on at least $4$ vertices or $\theta_0.$ To do so, it suffices to show that $\delta(F)>2.$ This holds since
$$
\delta(F)  = \delta(G|_Q) \ge \delta(G) +|Q| - m  = (2|Q| + 3\delta(G))/3+|Q|/3 -m \ge
$$
$$
(3m+2)/3 + 5/3 - m = 7/3 >2.
$$
II.2) Suppose that $F$ is not biconnected. Then there exists a cut vertex $v.$ We can argue as in Case I to show that 
$F|_{V(F)\backslash \{v\}}$ has exactly two connected components. Indeed, the degree of a vertex in $F|_{V(F)\backslash \{v\}}$ is at least $\delta(F)-1 = \delta(G|_Q)-1.$ Thus, for each connected component $K$ of $F|_{V(F)\backslash \{v\}},$ it holds that
$\delta(K)\ge \delta(G|_Q)-1\ge \delta(G)+|Q|-m -1,$ so the size of $K$ is at least $1+\delta(K) \ge \delta(G)+|Q|-m.$ As $3(\delta(G)+|Q|-m)> |Q|$ by the initial assumptions on $|Q|$ and $\delta(G),$ the claim that there are exactly two connected components follows.\\ 
Let $F_1'$ and $F_2'$ be the two connected components. We already know that they are both of minimum degree at least $\delta(G)+|Q|-m -1$ and of size at least $\delta(G)+|Q|-m.$ Now, pick any vertex $u\in V(F_2').$ Clearly, $u$ has no neighbors in $V(F_1'),$ so
$$
|V(F_1')|\le m - |N[u]| \le m-\delta(G)- 1.
$$
This means that
$$
\delta(F_1')\ge \delta(G)+|Q|-m -1 >(m-\delta(G)-1)/2\ge |V(F_1')|/2,
$$
again using the assumed inequality $2|Q|+3\delta(G)\ge 3m +2.$ By Lemma 2.6, the graph
$F_1'$ is Wilsonian. The same properties hold for $F_2'.$\\
Finally, suppose that $|N(v)\cap F_1'|\ge 2.$ Then clearly $F|_{V(F_1')\cup\{v\}}$ is also Wilsonian as it is still biconnected, not bipartite, not a cycle graph, and not isomorphic to $\theta_0.$ The same conclusion holds for $F_2'.$
\end{proof}\\

\subsection{Beginning of Main Proof}
We can assume without loss of generality that $\delta(X)\le \delta(Y)$ as $\FS(X,Y)$ and $\FS(Y,X)$ are isomorphic. This implies that $2\delta(X)+3\delta(Y)\ge 3n.$ 
Before we proceed to the main claims in the proof, we prove a simple bound on $\delta(X)+\delta(Y)$ that will be useful throughout:
$$
\delta(X)+\delta(Y) = 
(2\delta(X)+3\delta(Y))/3 + \delta(X)/3 >
3n/3 + n/6\ge n+1.
$$

We begin the main part of the proof using Lemma 2.11.  Fix a bijection
$\sigma: V(X)\longrightarrow V(Y)$ and vertices $u,v\in V(Y)$ such that 
$(\sigma^{-1}(u), \sigma^{-1}(v))\in E(X).$ We will prove that $u$ and $v$ are
$(X,Y)$-exchangeable from $\sigma.$ As $X$ is certainly connected (since $\delta(X)>|V(X)|/2$) and $u,v$ were chosen arbitrarily, 
it will then follow from Lemma 2.11 that $\FS(X, Y)$ is connected as desired.\\

Suppose, for the sake of contradiction, that $u$ and $v$ are not $(X,Y)$-exchangeable from $\sigma.$ This assumption will be maintained throughout all subsections until the end of Section 4 (i.e., in the proof of Theorem 1.4) and is an implicit condition of the intermediate lemmas (Lemma 4.2 - Lemma 4.5).\\

We proceed as follows. In Section 4.3, we will prove several claims on the sets of common neighbors of $u$ and $v$ in $X$ and $Y.$ These sets are natural to consider in light of the four techniques (especially number 4) outlined in Section 2.2.3 of the current paper. In Section 4.4, we will choose specific representatives of these sets and finish the proof via a combination of our four techniques.

\subsection{On the Common Neighbors of $u$ and $v$ in $X$ and $Y$}

Denote $u' = \sigma^{-1}(u), v' = \sigma^{-1}(v), A' = N(u')\cap N(v'), B = N(u)\cap N(v),A = \sigma(A'),$ and $B' = \sigma^{-1}(B).$\\

\noindent
\textbf{Lemma 4.2}:\textit{ The following inequalities hold:
$$
|B\cap \sigma(N(u'))|\le 1,$$
$$|B\cap \sigma(N(v'))|\le 1.
$$}
\begin{proof} We will only prove that $|B\cap \sigma(N(u'))|\le 1$ as the claims are symmetric.\\

Consider $\sigma(N[u']).$ Lemma 4.1 can be applied to $G = Y$ and  $Q = \sigma(N[u'])$ since $$|Q| = |N[u']| \ge \delta(X)+1>n/2+1\ge 4$$ and
$$2|Q|+3\delta(Y)\ge 2(\delta(X)+1)+3\delta(Y)\ge 3n+2.$$ From Lemma 4.1, it follows that $Y|_{\sigma(N[u'])}$ has at most two connected components. We consider several cases:\\

\noindent
\textbf{Case I:} The graph $Y|_{\sigma(N[u'])}$ consists of two connected components.\\
Let the components of $Y|_{\sigma(N[u'])}$ be $S_1$ and $S_2.$ By Lemma 4.1, we know that $S_1$ and $S_2$ are both Wilsonian. 
As $v'\in N(u'),$ it follows that $v$ is in one of the components, say, $S_1.$ We claim that $u\not \in V(S_1).$ Indeed, if $u$ is also in $S_1,$  then $\sigma^{-1}(S_1)$ contains a spanning star with center $u'.$ As $S_1$ is Wilsonian, Theorem 2.6 implies that  $u$ and $v$ are $(\sigma^{-1}(S_1),S_1)$- exchangeable from $\sigma|_{\sigma^{-1}(V(S_1))}.$ This means that they are also $(X,Y)$-exchangeable from $\sigma,$ which is a contradiction.\\
It follows that $u\in V(S_2).$ Thus, $N(u)\cap V(S_1) =\emptyset.$ Similarly,  $N(v)\cap V(S_2) = \emptyset$ as $v\in V(S_1).$ Combining these, we have 
$$B\cap \sigma(N(u')) \subseteq
(N(u)\cap N(v))\cap (V(S_1)\cup V(S_2))\subseteq 
(N(u)\cap V(S_1))\cup (N(v)\cap V(S_2))=\emptyset.
$$

\noindent
\textbf{Case II:} The graph $Y|_{\sigma(N[u'])}$ has a single connected component. Denote $S = Y|_{\sigma(N[u'])}$ for brevity. By Lemma 4.1, there are two cases:\\
II.1) Suppose that $S$ is almost-Wilsonian. Take an arbitrary $t'\in N(u').$ Then $$|N(t')\cap N(u')|\ge |N(t')|+|N(u')|-n\ge
2\delta(X)-n>0.
$$
Thus, $t'$ has a neighbor $s'\in N(u').$ Clearly 
$s'\neq u',$ so $X|_{N[u']}$ contains a spanning $\starplus_{|N[u']|}$ subgraph with center $u'$ and additional edge $(s', t').$ By Theorem 2.8, $u$ and $v$ are 
$(\sigma^{-1}(S),S)$-exchangeable from $\sigma|_{\sigma^{-1}(V(S))},$ so they are also $(X,Y)$-exchangeable from $\sigma,$ which is a contradiction.\\
II.2) Suppose that $S$ has a cut-vertex $x.$  By Lemma 4.1, $Y|_{V(S)\backslash \{x\}}$ has exactly two connected components $S_1$ and $S_2$ and both of them are Wilsonian.
We consider three cases:\\
II.2.1) Suppose that $x\not \in\{u,v\}.$ Without loss of generality, let $u$ be in the component $S_1.$ We claim that $v$ is not in $S_1.$ Indeed, $S_1$ is Wilsonian and $\sigma^{-1}(S_1)$ contains a spanning star with center $u'.$ Thus, if $v$ is also in $S_1,$ the vertices $u$ and $v$ will be $(\sigma^{-1}(S_1), S_1)$-exchangeable from 
$\sigma|_{V(\sigma^{-1}(S_1))}$ by Theorem 2.6. It follows that 
$u$ and $v$ are $(X,Y)$-exchangeable from $\sigma,$ which is a contradiction. Therefore, $v$ must be in $S_2.$ This clearly means that the only common neighbor of $u$ and $v$ in $V(S) = \sigma(N[u'])$ can be $x,$ so 
$$
1\ge |\sigma(N(u'))\cap (N(u)\cap N(v))| = 
|\sigma(N(u'))\cap B|,
$$
as desired.\\
II.2.2) Suppose that $x = v.$ It follows that $Y|_{V(S)\backslash \{v\}}$ has two connected components $S_1$ and $S_2.$ Without loss of generality, let $u\in S_1.$ We claim that  $|N(v)\cap V(S_1)|\le 1.$ Indeed, otherwise we know that 
$Y|_{V(S_1)\cup \{v\}}$ is Wilsonian from Lemma 4.1. As 
$X|_{\sigma^{-1}(V(S_1)\cup \{v\})}$ contains a spanning star subgraph with center $u',$ the vertices $u$ and $v$ are
$(X|_{\sigma^{-1}(V(S_1)\cup \{v\})}, Y|_{V(S_1)\cup \{v\}})$-exchangeable from $\sigma|_{\sigma^{-1}(V(S_1)\cup \{v\})}.$ Thus, they are also $(X,Y)$-exchangeable from $\sigma,$ which is a contradiction. Therefore, it must be the case that $|N(v)\cap V(S_1)|\le 1.$ This clearly implies that $u$ and $v$ have at most one common neighbor in
$V(S) = \sigma(N[u']).$ The inequality 
$1\ge |\sigma(N(u'))\cap B|$
is immediate.\\
II.2.3) Finally, in the case when $x = u,$ we can reason in an analogous way to Case II.2.2).
\end{proof}\\

\noindent
\textbf{Lemma 4.3}: \textit{ The following inequalities hold:
$$|B|\ge 2\delta(Y)+2-n,$$
$$|A'|\ge 2\delta(X) + 2\delta(Y)-2n.$$
}
\begin{proof}
First, note that $u$ and $v$ are not adjacent in $Y.$ Indeed, otherwise, the $(X,Y)$-friendly swap $uv$ exchanges $u$ and $v,$ which is a contradiction. Therefore:
$$
|B| = |N(u)\cap N(v)| = |N[u]\cap N[v]|\ge |N[u]|+|N[v]| -n \ge 2\delta(Y)+2-n,
$$
which proves the first inequality.\\
To prove the second inequality, we use Lemma 4.2. As $u$ and $v$ are not $(X,Y)$-exchangeable from $\sigma,$ we know that  $|B\cap \sigma(N(u'))|\le 1$ and 
$|B\cap \sigma(N(v'))|\le 1.
$ It follows that 
\newpage
\begin{equation*}
    \begin{split}
      |A'| = &|N(u')\cap N(v')| =  \\
& |\sigma(N(u'))\cap \sigma(N(v'))|\ge \\
& |(\sigma(N(u'))\cap(V(Y)\backslash B))\cap (\sigma(N(v'))\cap (V(Y)\backslash B))|\ge\\
& |\sigma(N(u'))\cap(V(Y)\backslash B)| + 
|\sigma(N(v'))\cap(V(Y)\backslash B)| - 
|V(Y)\backslash B|
\ge\\
& |\sigma(N(u'))|- 1 + |\sigma(N(v'))| -1 - |V(Y)\backslash B|\ge 
2\delta(X)-2 - (n -|B|)\ge \\
& 2\delta(X)-2 - n + 2\delta(Y)+2-n =\\
& 2\delta(X) + 2\delta(Y)-2n,
    \end{split}
\end{equation*}
as desired. \end{proof}\\

\subsection{Endgame}
We have the inequalities $|B\cap \sigma(N(u'))|\le 1$ and $|B\cap \sigma(N(v'))|\le 1$ from Lemma 4.2 and the inequality $|B|\ge 2\delta(Y)+2-n >2$ from Lemma 4.3. Combining these inequalities, we find a vertex $w\in B\backslash (\sigma(N(u'))\cup \sigma(N(v'))).$ Denote $w' = \sigma^{-1}(w).$\\

First, note that $w\not \in A = \sigma(N(u'))\cap \sigma(N(v'))$ as $w\not \in \sigma(N(u')).$ This, however, means that 
$$
|N(w)\cap A| = |N[w]\cap A|\ge|N[w]|+|A|-n \ge $$
$$
1+\delta(Y) + 2\delta(X)+2\delta(Y) - 2n  - n  = 1+2\delta(X)+3\delta(Y)-3n\ge 1,
$$
using the bound on $|A'| = |A|$ from Lemma 4.3. This means that $w$ has a neighbor in $A.$ Denote this neighbor by $x,$ and let $x' = \sigma^{-1}(x).$\\

Until the end of the proof, we will study the following two sets:
$$P = N(x')\cap N(w'),$$
$$R = N[w']\cap \sigma^{-1}(N[w]).$$ Note that $u'\not \in P$ and $v'\not \in P$ as $u'\not \in N(w'),v'\not \in N(w').$ Similarly, $u'\not \in R$ and  $v'\not \in R.$
Furthermore, observe that $x'$ and $w'$ are not adjacent, as otherwise the sequence of $(X,Y)$-friendly swaps:
$$
wx,wu,wv,wu,wx
$$
exchanges $u$ and $v.$ Thus, we find the following lower bound on the size of $P:$
$$
|P| = |N(x')\cap N(w')| = |N[x']\cap N[w']|\ge 
|N[x']|+|N[w']|-n \ge 2\delta(X)+2-n>2,
$$
so $|P|\ge 3.$
We also bound the size of $R:$
$$
|R| = |N[w']\cap \sigma^{-1}(N[w])|\ge |N[w']| + |\sigma^{-1}(N[w])| - n \ge$$
$$\delta(X)+1 +\delta(Y)+1-n = \delta(X) + \delta(Y)+2-n.
$$

\begin{figure}[!htb]
\begin{center}
   \includegraphics[width = .8\linewidth]{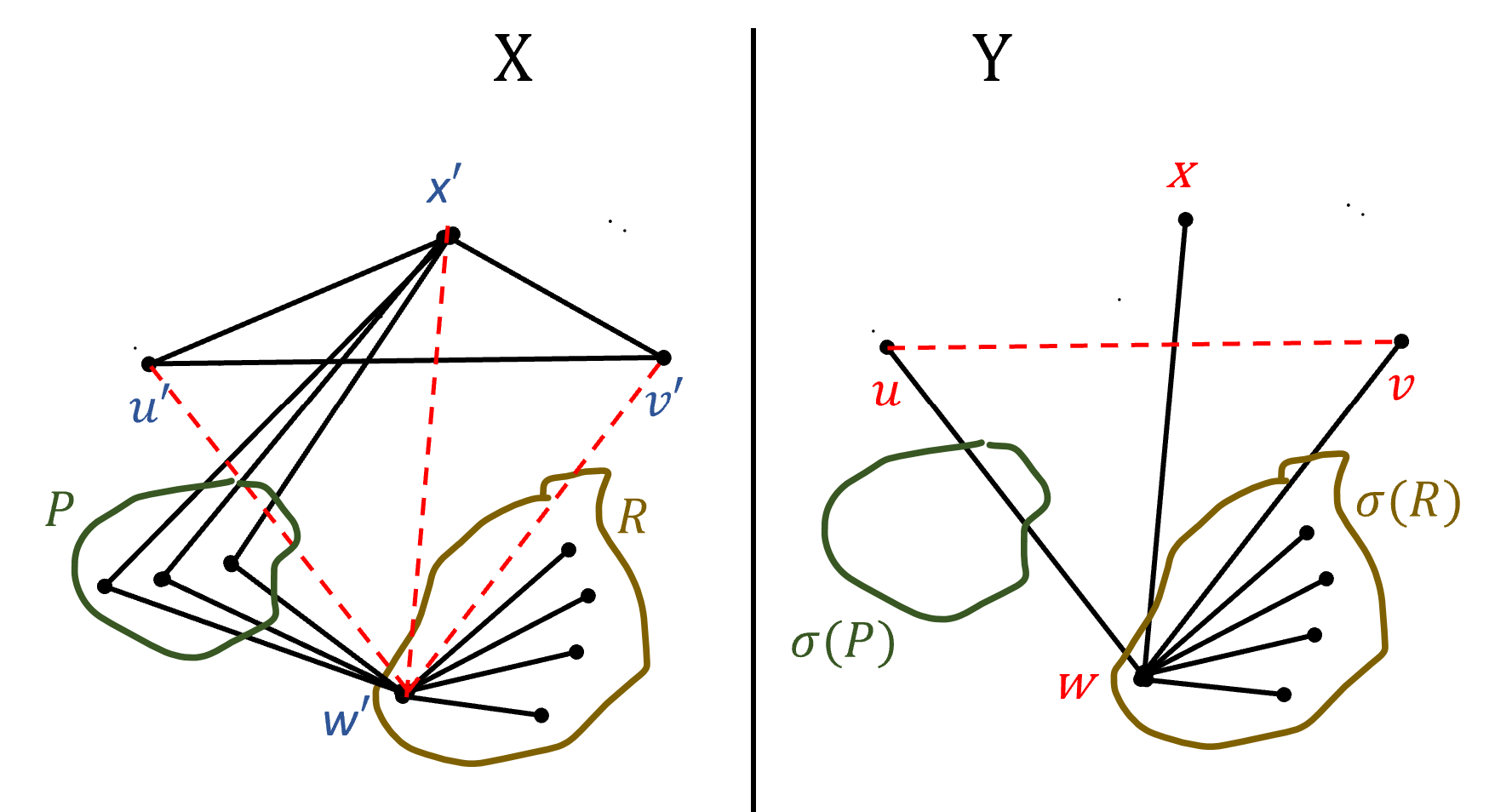}\caption{
   The vertices $w$ and $x$ and the sets $P$ and $R.$ In solid black are edges that exist and in dashed red are edges that do not exist.
   }\label{fig:Thm1.5.1}
\end{center}
\end{figure}

We illustrate the constructed vertices $x,w$ and the sets $P,R$ in \autoref{fig:Thm1.5.1}.\\

Next, we prove several short claims about $P$ and $R.$\\

\noindent
\textbf{Lemma 4.4}:\textit{ The sets $R$ and $N[x']$ are disjoint. In particular, $R$ and $P$ are disjoint.}\\
\begin{proof} Suppose, for the sake of contradiction, that there exists some $h' \in R\cap N[x'].$ Denote $h = \sigma(h').$ Note that $h'\not \in \{u',v'\}$ as $\{u',v'\}\cap R = \emptyset.$ Furthemore, $h'\in R\subseteq N[w']$ and $h'\in N[x']$ imply that $h'\not \in \{x',w'\}$ as $x'$ and $w'$ are not adjacent in $X.$ See \autoref{fig:Lemma4.4.1} for a diagram of the relevant vertices and edges.\\ 

\begin{figure}[!htb]
\begin{center}
   \includegraphics[width = .8\linewidth]{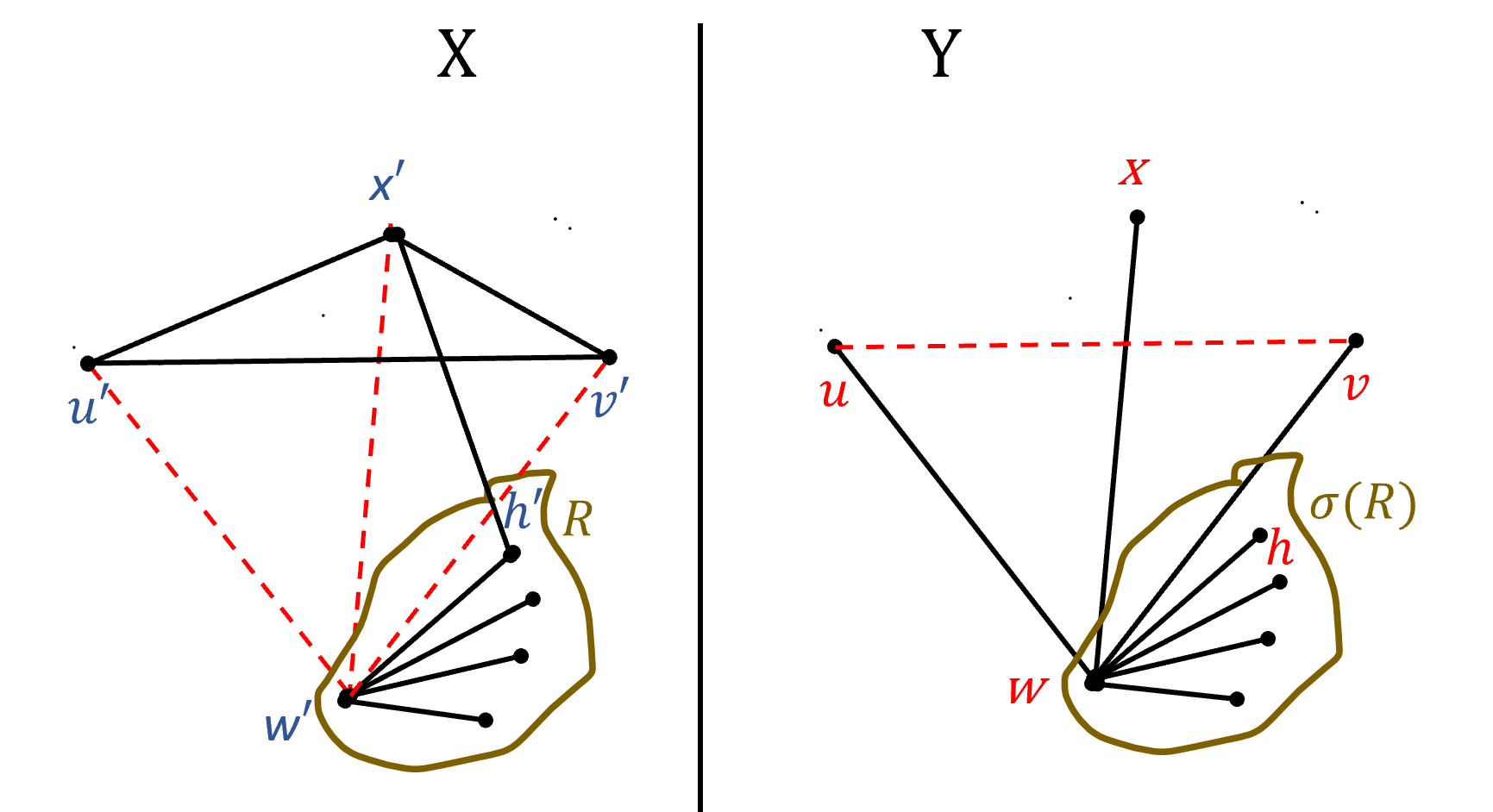}\caption{
   Diagram of $u,v,w,h,x$ and the relevant edges in $X$ and $Y.$ In solid black are edges that exist and in dashed red are edges that do not exist.
   }\label{fig:Lemma4.4.1}
\end{center}
\end{figure}

However, this means that the sequence of $(X,Y)$-friendly swaps
$$
wh, wx, wu, wv, wu, wx, wh
$$
exchanges $u$ and $v,$ which is a contradiction. Therefore, $R$ and $N[x']$ are indeed disjoint.
In particular, as $P\subseteq N[x'],$ the sets $P$ and $R$ are disjoint.
\end{proof}\\

\noindent
\textbf{Lemma 4.5}:\textit{
Every vertex in $ \sigma(P)$ has at most one neighbor in 
$\sigma(R).$}

\noindent
\begin{proof} Suppose, for the sake of contradiction, that there exists some $h\in \sigma(P)$ with at least two neighbors $q_1, q_2$ in $\sigma(R).$ Note that $q_1,q_2\in \sigma(R)\subseteq N[w].$ Denote 
$h' = \sigma^{-1}(h), q_1' = \sigma^{-1}(q_1),$ and $q_2 = \sigma^{-1}(q_2).$\\
First, note that neither of $q_1$ and $q_2$ is $w.$ Indeed, this would imply that $h\in N[w],$ and then $$h'\in P\cap \sigma^{-1}(N[w]) = N(x')\cap N(w')\cap \sigma^{-1}(N[w])\subseteq  N(w')\cap \sigma^{-1}(N[w])\subseteq R,
$$
whereas $P$ and $R$ are disjoint by Lemma 4.4. See \autoref{fig:Lemma4.5.1} below for a diagram of the relevant vertices and edges. \\\

\begin{figure}[!htb]
\begin{center}
   \includegraphics[width = .75\linewidth]{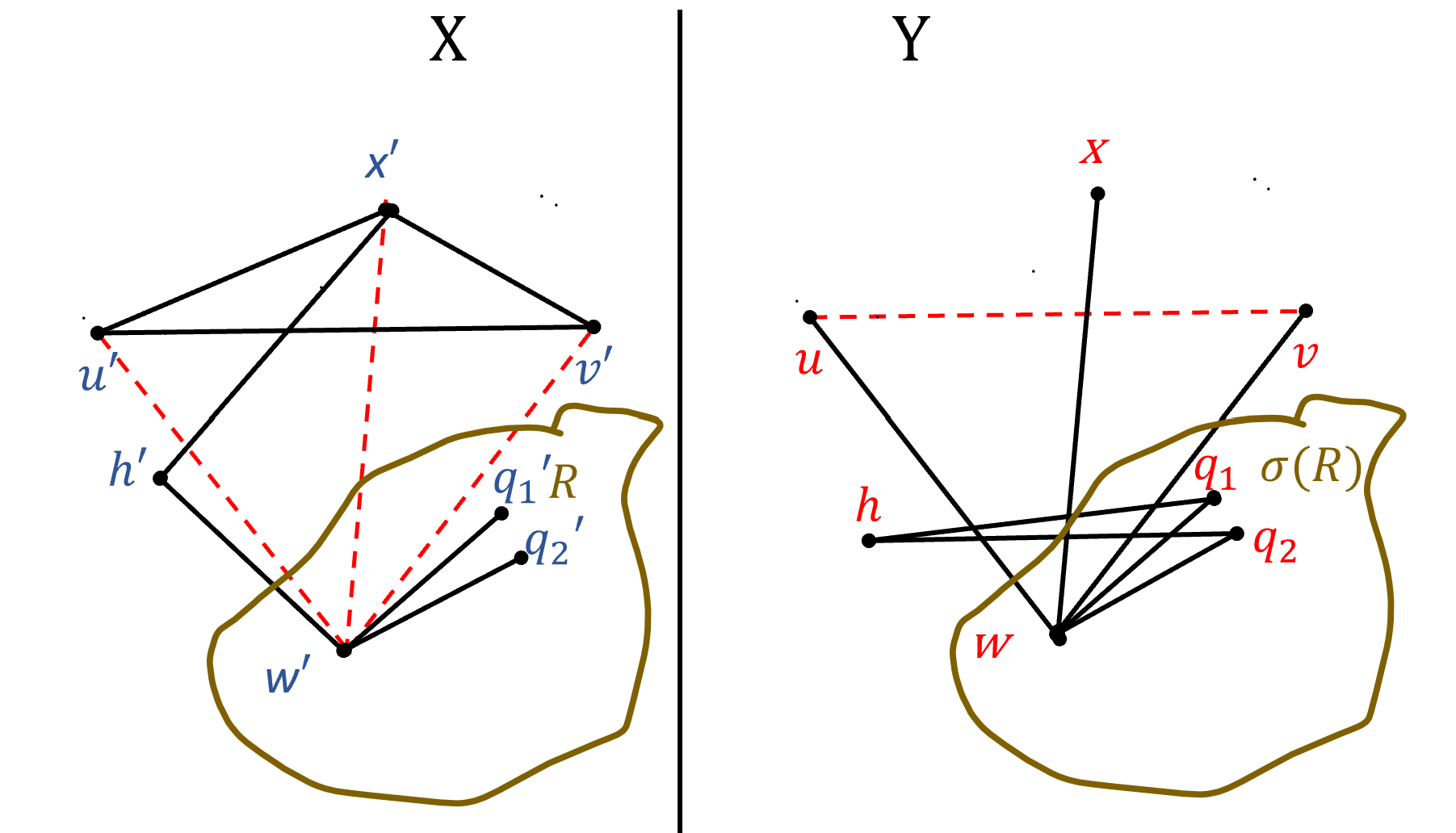}\caption{
   Diagram of $u,v,w,h,x,q_1, q_2$ and the relevant edges in $X$ and $Y.$ Again, in solid black are edges that exist and in dashed red are edges that do not exist.
   }\label{fig:Lemma4.5.1}
\end{center}
\end{figure}

Thus, $q_1, q_2, w$ are different vertices in $\sigma(R).$ Finally, the sequence of $(X,Y)$-friendly swaps
$$wq_1, hq_1, hq_2,wq_2, wq_1,wx, wu, wv, wu, wx, wq_1, wq_2, hq_2, hq_1, wq_1
$$
exchanges $u$ and $v,$ which is a contradiction. \end{proof}
\\

Now, we will finish the proof using Lemmas 4.4 and 4.5. First, note that Lemma 4.1 can be applied to $G = Y$ and $Q = \sigma(N[x'])$ since
$$
|Q|= |\sigma(N[x'])|\ge 1 + \delta(X)>1+n/2\ge 4,
$$
and
$$2|\sigma(N[x'])|+3\delta(Y)\ge
2(\delta(X)+1)+3\delta(Y)\ge 3n+2.$$
Thus, by Lemma 4.1, the graph  $Y|_{\sigma(N[x'])}$ consists of at most two connected components. We consider two cases based on the number of connected components. \\

\noindent
\textbf{Case I:} Suppose that $Y|_{\sigma(N[x'])}$ has two connected components, $F_1$ and $F_2.$
Take any vertex $q\in \sigma(P)\subseteq \sigma(N[x']).$ Without loss of generality, $q\in F_1.$ As $F_2$ and $F_1$ are distinct connected components,
we know that $N(q)\cap V(F_2) = \emptyset.$ 
Furthermore, as $V(F_2)\subseteq\sigma(N[x']),$ Lemma 4.4 implies that $\sigma(R)\cap V(F_2) = \emptyset.$ 
At the same time, by Lemma 4.1, $$|V(F_2)|\ge \delta(Y)+1 + |\sigma(N[x'])|-n\ge \delta(Y)+\delta(X)+2-n.$$ Alos, using the fact 
$|N(q)\cap \sigma(R)|\le 1$ from Lemma 4.5, we conclude that $q$ has at most 1 neighbor in the entire set $V(F_2)\cup \sigma(R).$ However, the degree of $q$ in $Y$ is too large for this to be possible. To formalize this idea, we compute
\begin{equation*}
    \begin{split}
        n  \ge & |N(q)\cup\sigma(R)\cup V(F_2)|=\\
        & |N(q)\cup \sigma(R)| + 
|V(F_2)|\ge\\
& |N(q)| + |\sigma(R)|-1 + (\delta(Y)+\delta(X)+2-n)\ge\\
& \delta(Y)+(\delta(X)+\delta(Y)+2-n)-1 
+
(\delta(Y)+\delta(X)+2-n)\ge\\
& 3\delta(Y)+2\delta(X)-2n+3\ge n+3,
    \end{split}
\end{equation*}
which leads to the 
contradiction $n\ge n+3.$\\

\noindent
\textbf{Case II:} Suppose that $Y|_{\sigma(N[x'])}$ has a single connected component $F.$ We consider two cases:\\
II.1) Suppose that $F$ is almost-Wilsonian. Take an arbitrary $t'\in N(x').$ Then
$$
|N(t')\cap N(x')|\ge 
|N(t')|+|N(x')|-n \ge 2\delta(X)-n>0,
$$
so $t'$ has a neighbor $s'$ in $N(x').$ Clearly, $t'\neq x'.$ 
Therefore, 
$X|_{\sigma^{-1}(V(F))}$ contains a spanning $\starplus_{|N[x']|}$ graph with center $x'$ and additional edge $(s',t').$ As $u,v\in V(F) =\sigma(N[x']),$ Theorem 2.8 implies that the vertices $u$ and $v$ are  $(X|_{\sigma^{-1}(V(F))}, F)$-exchangeable from 
$\sigma|_{\sigma^{-1}(V(F))}.$ Thus, they are $(X,Y)$-exchangeable from $\sigma,$ which is a contradiction.\\
II.2) Suppose that $F$ is not almost-Wilsonian. Then it has a cut vertex $y.$ Furthermore, $F|_{V(F)\backslash\{y\}}$ has two connected components $F_1$ and $F_2,$ each of size at least 
$
\delta(Y)+|N[x']|-n\ge \delta(Y)+\delta(X)+1-n.
$
As $|P|> 2,$
there exists a vertex $q\in \sigma(P)$ different from $y.$ Without loss of generality, $q\in F_1.$ This implies that
$N(q)\cap V(F_2)=\emptyset.$ As
$|N(q)\cap \sigma(R)|\le 1$ by Lemma 4.5
and $V(F_2)\subseteq \sigma(N[x'])$ is disjoint from $\sigma(R)$ by Lemma 4.4, we argue as in Case I:
\begin{equation*}
    \begin{split}
      n \ge & |N(q)\cup \sigma(R)\cup V(F_2)| = \\
      & |N(q)\cup \sigma(R)|+|V(F_2)|\ge\\
      & (|N(q)| + |\sigma(R)| - 1) + 
\delta(X)+\delta(Y)+1-n
\ge\\
    & \delta(Y) + (\delta(X)+\delta(Y)+2 - n) - 1+ 
\delta(X)+\delta(Y)+1-n\ge\\
& 3\delta(Y)+2\delta(X)+2-2n\ge n+2,\\
    \end{split}
\end{equation*}
which is our final contradiction.

\section{Lower Bounds for Arbitrary Graphs: Proposition 1.6}
\noindent
\textbf{Proposition 1.6}:\textit{ Suppose that $n\ge k\ge 5$ are integers. Then there exist connected graphs $X$ and $Y$ such that $\delta(X)\ge \frac{3n}{k}-4,\delta(Y)\ge \frac{(k-2)n}{k}-3$ and $\FS(X,Y)$ is disconnected.}\\
\begin{proof} Let $n = kt + r,$ where $r$ and $t$ are integers and $0\le r \le k-1.$ \\

We first construct $X.$ Split $V(X)$ into $k$ groups  $A_1, A_2, \ldots , A_k$ of almost equal size, that is, $|A_i| = t+1$ for $1 \le i \le r$ and $|A_i| = t$ for $i >r.$ If $u$ and $v$ are distinct vertices such that $u \in A_i$ and $v\in A_j,$ we connect $u$ and $v$ if and only if $i - j \in \{-1,0,1\}\pmod{k}.$ It is easy to check that the minimum degree of $X$ is at least $3t-1 \ge \frac{3n}{k}-4.$\\

Now, we construct $Y.$ Split $V(Y)$ into $k$ groups  $B_1, B_2, \ldots , B_k,$ such that $|B_i| = t+1$ for $1 \le i \le r$ and $|B_i| = t$ for $i >r.$ If $u$ and $v$ are distinct vertices such that $u \in B_i$ and $v\in B_j,$ we connect $u$ and $v$ if and only if $i - j \not \in \{-1,1\}\pmod{k}.$ Again, one can check that the minimum degree of $Y$ is at least $n - (2t-3)\ge \frac{(k-2)n}{k} - 3.$ The graphs $X$ and $Y$ are illustrated in \autoref{fig:LowerBoundFS}.\\

\begin{figure}[!htb]
\begin{center}
   \includegraphics[width = .75\linewidth]{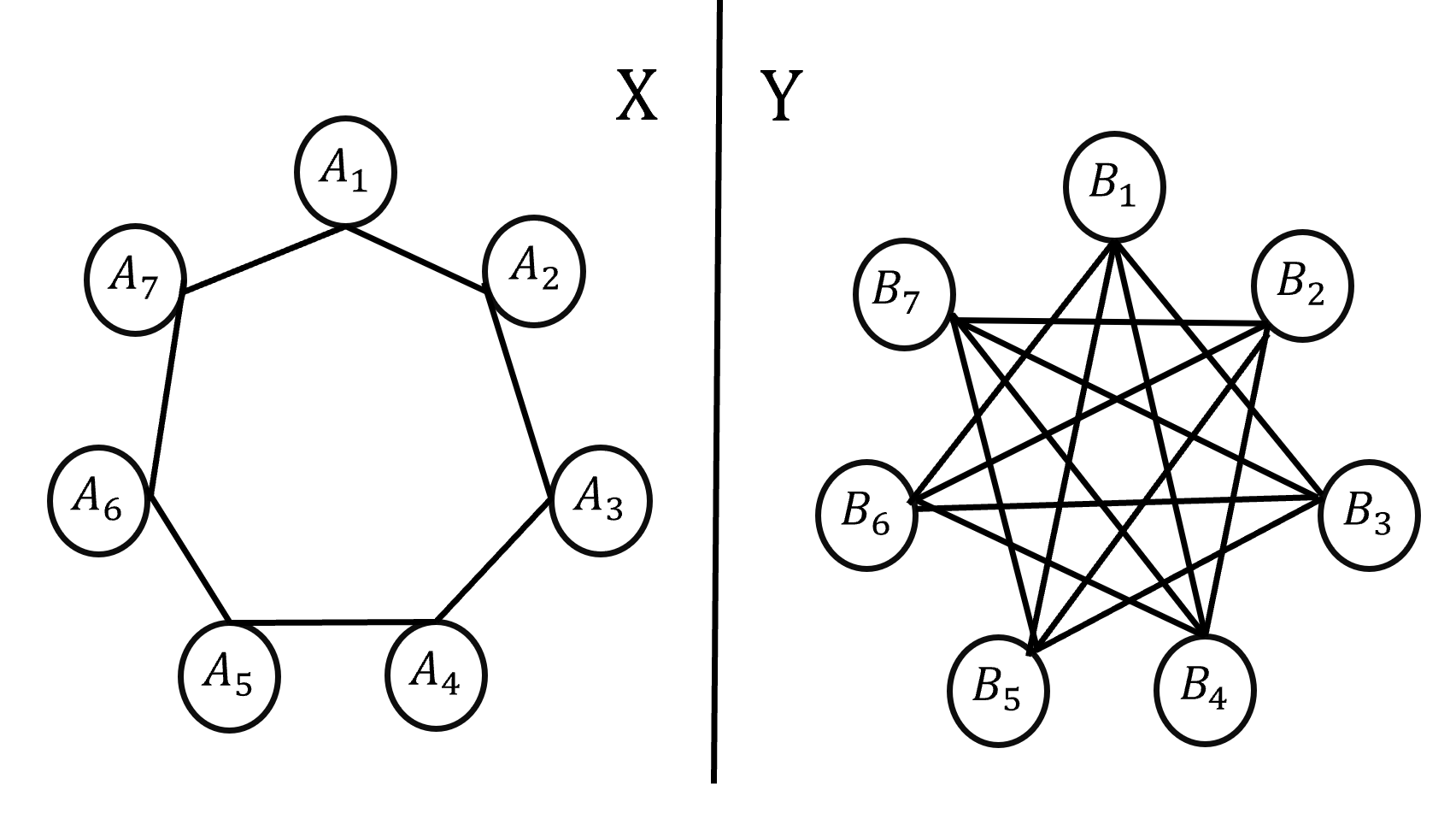}\caption{
   Illustration of the graphs $X$ and $Y$ from Proposition 1.6 in the case $k = 7.$ The induced subgraphs on each $A_i$ and $B_i$ are cliques. Two sets $A_i$ and $A_j$ (respectively $B_i$ and $B_j$) are connected on the diagram if all edges between $A_i$ and $A_j$ (respectively $B_i$ and $B_j$) exist. There are no other edges.
   }\label{fig:LowerBoundFS}
\end{center}
\end{figure}

To show that $\FS(X,Y)$ is disconnected, consider a bijection $\sigma: V(X)\longrightarrow V(Y)$ satisfying $\sigma(A_i) = B_i$ for each $i.$ It is clear that $uv$ can be an $(X,Y)$-friendly swap only if $u$ and $v$ are in the same $B_j.$ It easily follows that any other bijection $\tau$ in the connected component of $\sigma$ in $\FS(X,Y)$ will also satisfy that $\tau(A_i) = B_i$ for each $i.$ Clearly, not all bijections from $V(X)$ to $V(Y)$ satisfy this property, so $\FS(X,Y)$ is disconnected.\end{proof}

\section{Upper Bounds for Bipartite Graphs: Theorem 1.10}

\noindent
\textbf{Theorem 1.10}: \textit{Let $r\ge 2,$ and let $X$ and $Y$ be edge-subgraphs of $K_{r,r}$ such that $$\delta(X)+\delta(Y)\ge 3r/2+1.$$ Then $\FS(X,Y)$ has exactly two connected components}.\\

To prove Theorem 1.10, it is sufficient to establish the following proposition.\\

\noindent
\textbf{Proposition 6.1:}
\textit{Let $r\ge 2,$ and let $X$ and $Y$ be edge-subgraphs of 
$K_{r,r}$ such that $\delta(X)+\delta(Y)\ge 3r/2 + 1.$ Let $\sigma: V(X)\longrightarrow V(Y)$ be an arbitrary bijection. If $u,v\in V(Y)$ are in different partite sets of $Y$ and are such that $\{\sigma^{-1}(u), \sigma^{-1}(v)\}\in E(X),$ then $u$ and $v$ are $(X,Y)$-exchangeable from $\sigma.$} \\

Before we prove the proposition, let us see how it implies Theorem 1.10. Suppose $X$ and $Y$ are edge-subgraphs of $K_{r,r}$ with $\delta(X)+\delta(Y)\ge 3r/2 + 1.$ We want to show that $\FS(X,Y)$ has two connected components. Proposition 6.1 tells us that the hypothesis of Proposition 2.10 is satisfied with $\widetilde{Y} = K_{r,r},$ so it follows that $\FS(X,Y)$ and $\FS(X,K_{r,r})$ have the same number of connected components. Recall that
$\FS(X,K_{r,r})\cong \FS(K_{r,r}, X)$ by Proposition 2.1, and denote $X' = K_{r,r}, Y' = X.$ We can apply Proposition 6.1 to $X'$ and $Y'$ since $\delta(X') + \delta(Y') = r +\delta(X)\ge \delta(Y)+\delta(X)\ge 3r/2 + 1.$ Therefore, the number of connected components of $\FS(X',Y')$ is the same as the number of connected components of 
$\FS(X', K_{r,r}) \cong \FS(K_{r,r}, K_{r,r}),$ which we know is two by Proposition 2.3. Thus, $\FS(X,Y)$ has two connected components whenever $\delta(X)+\delta(Y)\ge 3r/2 + 1.$\\

In the rest of Section 6, we will prove Proposition 6.1. We split the proof in multiple subsections for clarity. Throughout all of them, all assumptions on $X,Y,u,v,$ and $\sigma$ given in the statement of Proposition 6.1 are maintained.

\subsection{Beginning of the Proof}

Let $\{A_X,B_X\}$ and $\{A_Y, B_Y\}$ be the bipartitions of $X$ and $Y,$ respectively. Without loss of generality, we may assume that $u \in A_Y$ and $v\in B_Y.$ Let $u' = \sigma^{-1}(u)$ and $v' = \sigma^{-1}(v).$ Our goal is to show that $\sigma$ and $(u,v)\circ\sigma$ are in the same connected component of $\FS(X,Y).$ We may assume that the partite set of $X$ containing $u'$ contains at least $(r-1)/2$ elements of $\sigma^{-1}(A_Y\backslash \{u\})$ since otherwise we can switch the roles of $\sigma$ and $(u,v)\circ \sigma.$ Without loss of generality, we may assume that $u'\in A_X$ and $v'\in B_X.$ Thus, $|(A_X\backslash\{u'\})\cap \sigma^{-1}(A_Y\backslash \{u\})|\ge (r-1)/2$ holds.

\subsection{Main Lemma}
\noindent
\textbf{Lemma 6.2:} \textit{There exists a sequence $\Sigma$ of $(X,Y)$-friendly swaps which does not involve $u$ or $v$ such that $\Sigma$ transforms $\sigma$ into a bijection $\mu$ satisfying $\mu(A_X) = A_Y$ and $\mu(B_X) = B_Y.$}\\ 
\begin{proof}
We first introduce the edge-subgraphs of $X$ and $Y$ given by $\widetilde{X} = X|_{V(X)\backslash \{u',v'\}}$ and $\widetilde{Y} = Y|_{V(Y)\backslash \{u,v\}}.$ Note that $\sigma(V(\widetilde{X})) = V(\widetilde{Y}).$ 
Let the partite sets of $\widetilde{X}$ be $A_{\widetilde{X}}$ and $B_{\widetilde{X}},$ each of size $r-1,$ and the partite sets of $\widetilde{Y}$ be $A_{\widetilde{Y}}$ and $B_{\widetilde{Y}},$ again each of size $r-1.$ Observe that $\delta(\widetilde{X})\ge \delta(X)-1$ and $\delta(\widetilde{Y}) \ge \delta(Y)-1.$ In particular, this means that $\delta(\widetilde{X})+\delta(\widetilde{Y})\ge 3r/2-1.$
\\

In order to show that there exists a sequence $\Sigma$ of $(X,Y)$-friendly swaps which does not involve $u$ or $v$ and transforms $\sigma$ into a bijection $\mu$ satisfying $\mu(A_X) = A_Y,$ it is sufficient to show that there exists a bijection $\widetilde{\mu}$ in the same connected component of $\FS(\widetilde{X},\widetilde{Y})$ as $\sigma|_{V(\widetilde{X})}$ satisfying $\widetilde{\mu}(A_{\widetilde{X}}) = A_{\widetilde{Y}}.$\\ 

Consider a bijection $\widetilde{\mu}$ in the same connected component of $\FS(\widetilde{X},\widetilde{Y})$ as $\sigma|_{V(\widetilde{X})}$ for which the quantity $|\widetilde{\mu}(A_{\widetilde{X}})\cap A_{\widetilde{Y}}|$ is maximal. Let $s = |\widetilde{\mu}(A_{\widetilde{X}})\cap A_{\widetilde{Y}}|.$ We want to show that $s = r-1.$ Suppose, for the sake of contradiction, that $s<r-1.$
Note that $s\ge (r-1)/2$ holds because 
$$
|\sigma|_{V(\widetilde{X})}(A_{\widetilde{X}})\cap A_{\widetilde{Y}}| = 
|(A_X\backslash\{u'\})\cap \sigma^{-1}(A_Y\backslash \{u\})|\ge (r-1)/2
$$
as described above.\\

Now, observe that $|\widetilde{\mu}(A_{\widetilde{X}})\cap A_{\widetilde{Y}}| = s$ implies that
$$
|\widetilde{\mu}(A_{\widetilde{X}})\cap B_{\widetilde{Y}}| = r-1-s, \; \; |\widetilde{\mu}(B_{\widetilde{X}})\cap A_{\widetilde{Y}}| = r-1-s, \;\;\text{and  }\;\; |\widetilde{\mu}(B_{\widetilde{X}})\cap B_{\widetilde{Y}}| = s.
$$
Since $s<r-1,$ the sets $\widetilde{\mu}(A_{\widetilde{X}})\cap B_{\widetilde{Y}}$ and $\widetilde{\mu}(B_{\widetilde{X}})\cap A_{\widetilde{Y}}$ are both non-empty.
We consider two cases based on the edges between those sets.\\
\textbf{Case I:} Suppose that there exist two vertices $a\in \widetilde{\mu}(A_{\widetilde{X}})\cap B_{\widetilde{Y}}$ and $b\in \widetilde{\mu}(B_{\widetilde{X}})\cap A_{\widetilde{Y}}$
such that $(a,b)\in E(Y)$ and
$(\widetilde{\mu}^{-1}(a), \widetilde{\mu}^{-1}(b))\in E(X).$ Then $ab$ is an $(\widetilde{X}, \widetilde{Y})$-friendly swap from $\widetilde{\mu}.$ Therefore, 
$\widetilde{\nu} = (a,b)\circ \widetilde{\mu}$ is in the same connected component of $\FS(\widetilde{X},\widetilde{Y})$ as $\widetilde{\mu}.$ This contradicts the maximality of $|\widetilde{\mu}(A_X)\cap A_Y|$ since $$\widetilde{\nu}(A_X\backslash \{\widetilde{\mu}^{-1}(a)\}) = \widetilde{\mu}(A_X\backslash \{\widetilde{\mu}^{-1}(a)\}),$$ $$\widetilde{\nu}(\widetilde{\mu}^{-1}(a)) = (a,b)[\widetilde{\mu}(\widetilde{\mu}^{-1}(a))] = b \in A_{\widetilde{Y}},\; \text{ and  }\;\;
\widetilde{\mu}(\widetilde{\mu}^{-1}(a)) = a\in B_Y,$$ which means that
$$
|\widetilde{\nu}(A_{\widetilde{X}})\cap A_{\widetilde{Y}}| = |\widetilde{\mu}(A_{\widetilde{X}})\cap A_{\widetilde{Y}}|+1.
$$
\textbf{Case II:} Suppose that no two vertices $p\in \widetilde{\mu}(A_{\widetilde{X}})\cap B_{\widetilde{Y}}$ and $q\in \widetilde{\mu}(B_{\widetilde{X}})\cap A_{\widetilde{Y}}$ simultaneously satisfy that 
$(p,q)\in E(Y)$ and $(\widetilde{\mu}^{-1}(p), \widetilde{\mu}^{-1}(q))\in E(X).$ Let the number of edges between $\widetilde{\mu}(A_{\widetilde{X}})\cap B_{\widetilde{Y}}$ and $\widetilde{\mu}(B_{\widetilde{X}})\cap A_{\widetilde{Y}}$ be $m.$ 
Since $|\widetilde{\mu}(A_{\widetilde{X}})\cap B_{\widetilde{Y}}| =|\widetilde{\mu}(B_{\widetilde{X}})\cap A_{\widetilde{Y}}|= r-1-s,$ 
the assumption in Case II implies that there are at most $(r-1-s)^2 - m$ edges between $\widetilde{\mu}^{-1}(\widetilde{\mu}(A_{\widetilde{X}})\cap B_{\widetilde{Y}}) = A_{\widetilde{X}}\cap \widetilde{\mu}^{-1}(B_{\widetilde{Y}})$ and 
$
\widetilde{\mu}^{-1}(\widetilde{\mu}(B_{\widetilde{X}})\cap A_{\widetilde{Y}}) = 
B_{\widetilde{X}}\cap \widetilde{\mu}^{-1}(A_{\widetilde{Y}}).
$\\

Let $a$ be a vertex in $\widetilde{\mu}(A_{\widetilde{X}})\cap B_{\widetilde{Y}}$ with the fewest neighbors in $\widetilde{\mu}(B_{\widetilde{X}})\cap A_{\widetilde{Y}}.$ By the pigeonhole principle, 
$$
|N(a)\cap (\widetilde{\mu}(B_{\widetilde{X}})\cap A_{\widetilde{Y}})|\le \frac{m}{r-1-s}.
$$
We will denote $t = m/(r-1-s)$ for brevity. So $|N(a)\cap (\widetilde{\mu}(B_{\widetilde{X}})\cap A_{\widetilde{Y}})|\le t.$
Similarly, let $b'$ be a vertex in 
$B_{\widetilde{X}}\cap \widetilde{\mu}^{-1}(A_{\widetilde{Y}})$ with the fewest neighbors in $A_{\widetilde{X}}\cap \widetilde{\mu}^{-1}(B_{\widetilde{Y}}).$ By the pigeonhole principle, 
$$
|N(b')\cap (A_{\widetilde{X}}\cap \widetilde{\mu}^{-1}(B_{\widetilde{Y}}))|\le \frac{(r-1-s)^2 - m}{r-1-s} = (r-1-s) - \frac{m}{r-1-s} = r-1-s-t.
$$

Denote $a' = \widetilde{\mu}^{-1}(a)$ and 
$b = \widetilde{\mu}(b').$\\

Now, we will show that there exists a vertex $c'\in B_{\widetilde{X}}\cap \widetilde{\mu}^{-1}(B_{\widetilde{Y}})$ such that $(a',c')\in E(X)$ and $(b,\widetilde{\mu}(c')) \in E(Y).$ 
The existence of such a $c'$ will clearly follow from the inequality 
$$|N(a')\cap (B_{\widetilde{X}}\cap \widetilde{\mu}^{-1}(B_{\widetilde{Y}}))|+
|N(b)\cap (\widetilde{\mu}(B_{\widetilde{X}})\cap B_{\widetilde{Y}})|> |B_{\widetilde{X}}\cap \widetilde{\mu}^{-1}(B_{\widetilde{Y}})|.$$
To prove this inequality, we bound
$$
|N(a')\cap (B_{\widetilde{X}}\cap \widetilde{\mu}^{-1}(B_{\widetilde{Y}}))| = 
|N(a')\cap B_{\widetilde{X}}| - 
|N(a')\cap (B_{\widetilde{X}}\cap \widetilde{\mu}^{-1}(A_{\widetilde{Y}})|\ge
$$
$$
\delta(\widetilde{X}) - |B_{\widetilde{X}}\cap \widetilde{\mu}^{-1}(A_{\widetilde{Y}})| = 
\delta(\widetilde{X}) - (r-1-s).
$$
Similarly, 
$$
|N(b')\cap (\widetilde{\mu}(B_{\widetilde{X}})\cap B_{\widetilde{Y}})|\ge
\delta(\widetilde{Y}) - (r-1-s).
$$
As $|B_{\widetilde{X}}\cap \widetilde{\mu}^{-1}(B_{\widetilde{Y}})| = s,$ we only need to show that 
$$
\delta(\widetilde{X}) - (r-1-s)+\delta(\widetilde{Y}) - (r-1-s)> s.
$$
This is equivalent to showing that
$
\delta(\widetilde{X}) + \delta(\widetilde{Y}) + s+2> 2r.
$
The last inequality follows from the inequalities 
$s\ge (r-1)/2$ and $\delta(\widetilde{X})+\delta(\widetilde{Y})\ge 3r/2-1,$ which we already proved. Thus, such a vertex $c'$ exists. Denote $c = \widetilde{\mu}(c').$\\

Next, we will show that there exists a vertex $d\in \widetilde{\mu}(A_{\widetilde{X}})\cap A_{\widetilde{Y}}$ which satisfies the following four properties:
\vspace*{-0.8em}
\begin{itemize}
    \setlength\itemsep{-0.4em}
    \item $(c,d)\in E(Y),$
    \item $(a,d)\in E(Y),$
    \item $(c', \widetilde{\mu}^{-1}(d))\in E(X),$ 
    \item $(b',\widetilde{\mu}^{-1}(d))\in E(X).$
\end{itemize}

To prove the existence of such a vertex, it is clearly sufficient to show that 
\begin{equation}
    \tag{6.1}
|N(a)\cap N(c)\cap (\widetilde{\mu}(A_{\widetilde{X}})\cap A_{\widetilde{Y}})|+
|N(b')\cap N(c')\cap (A_{\widetilde{X}}\cap \widetilde{\mu}^{-1}(A_{\widetilde{Y}}))| > 
|\widetilde{\mu}(A_{\widetilde{X}})\cap A_{\widetilde{Y}}|.
\end{equation}
We first bound the number of common neighbors of $a$ and $c$ in $\widetilde{\mu}(A_{\widetilde{X}})\cap A_{\widetilde{Y}}.$ Note that
$$
|(N(a)\cap A_{\widetilde{Y}})\cap 
(N(c)\cap A_{\widetilde{Y}})| \ge 
2\delta(\widetilde{Y}) - |A_{\widetilde{Y}}| = 
2\delta(\widetilde{Y}) - (r-1).
$$
Since $a$ has at most $t$ neighbors in $\widetilde{\mu}(B_{\widetilde{X}})\cap A_{\widetilde{Y}},$ it follows that $a$ and $c$ have at least $2\delta(\widetilde{Y}) - (r-1) - t$ common neighbors in 
$A_{\widetilde{Y}}\backslash (\widetilde{\mu}(B_{\widetilde{X}})\cap A_{\widetilde{Y}}) = \widetilde{\mu}(A_{\widetilde{X}})\cap A_{\widetilde{Y}}.$\\

In the same way, we find that $b'$ and $c'$ have at least $2\delta(\widetilde{X}) - (r-1) - (r-1-s-t)$ common neighbors in $A_{\widetilde{X}}\cap \widetilde{\mu}^{-1}(A_{\widetilde{Y}}).$ \\

Since $|\widetilde{\mu}(A_{\widetilde{X}})\cap A_{\widetilde{Y}}| = s,$ to prove the desired inequality (6.1),
it is sufficient to show that
$$
2\delta(\widetilde{Y}) - (r-1) - t + 
2\delta(\widetilde{X}) - (r-1) - (r-1-s-t)> s.
$$
The latter is equivalent to
$$
2\delta(\widetilde{Y}) + 2\delta(\widetilde{X}) > 
3r-3,
$$
which follows from the inequality $\delta(\widetilde{X})+\delta(\widetilde{Y})\ge 3r/2-1.$ The existence of such a vertex $d$ follows. Let $d' = \widetilde{\mu}^{-1}(d).$\\

So far, we have found vertices $a\in \widetilde{\mu}(A_{\widetilde{X}})\cap B_{\widetilde{Y}}, b \in \widetilde{\mu}(B_{\widetilde{X}})\cap A_{\widetilde{Y}},c\in \widetilde{\mu}(B_{\widetilde{X}})\cap B_{\widetilde{Y}},$ and 
$d\in \widetilde{\mu}(A_{\widetilde{X}})\cap A_{\widetilde{Y}},$ such that the following edges between them exist:
\vspace*{-0.8em}
\begin{itemize}
    \setlength\itemsep{-0.4em}
    \item $(a,d),(b,c),(c,d)\in E(Y),$
    \item $(a',c'),
    (b',d'),
    (c',d')
     \in E(X).$
\end{itemize}
See \autoref{fig:Bipartite_d_rr_full} for a visualization of these edges and vertices.\\

\begin{figure}[!htb]
\begin{center}
   \includegraphics[width =.8\linewidth]{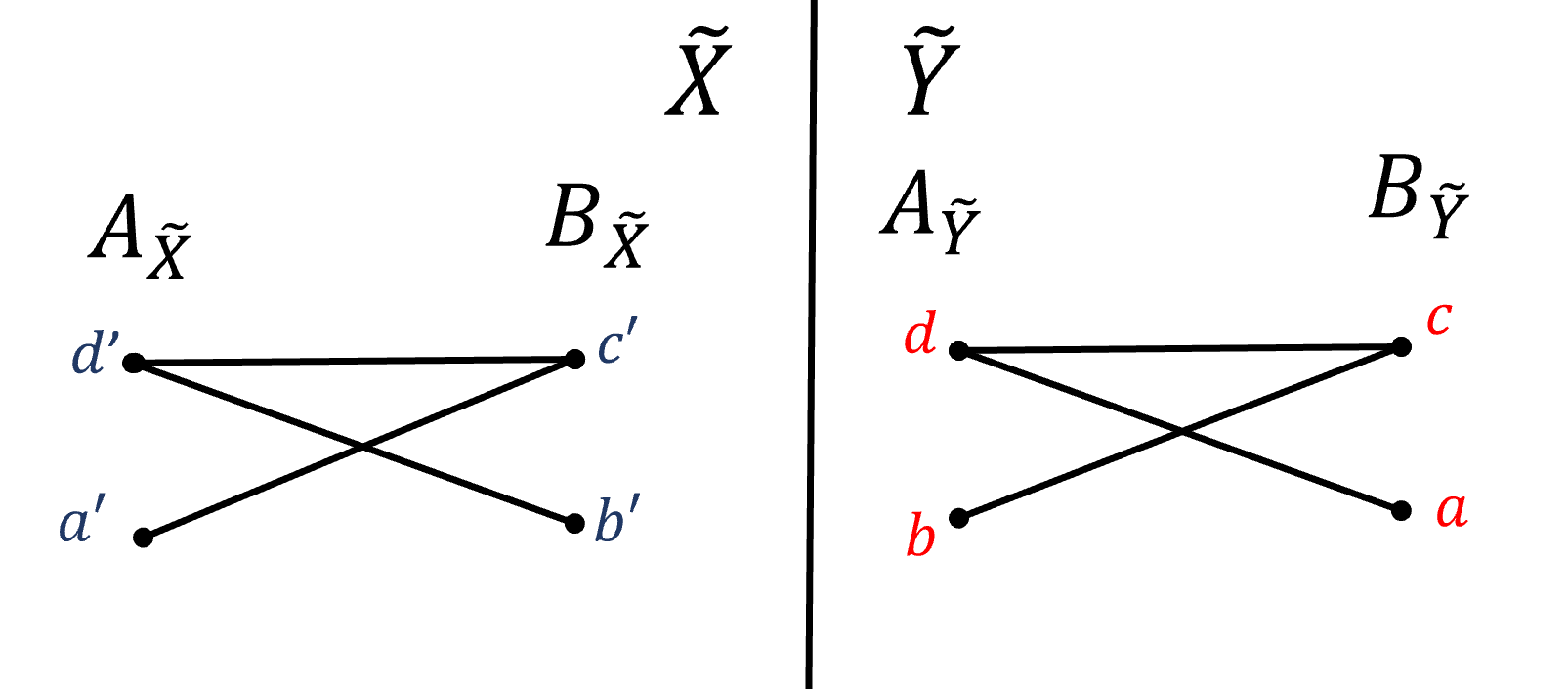}\caption{
   Relevant edges and vertices to Case II in the proof of Lemma 6.2.
   }\label{fig:Bipartite_d_rr_full}
\end{center}
\end{figure}

To finish the argument in Case II, we consider the sequence of $(\widetilde{X}, \widetilde{Y})$-friendly swaps
$$ 
cd, ad, bc
$$
applied to $\widetilde{\mu}.$ Let the resulting bijection be $\widetilde{\nu}.$ Observe that $\widetilde{\mu}$ and $\widetilde{\nu}$ satisfy
\begin{equation*}
    \begin{split}
    & \widetilde{\mu}^{-1}(a)\in A_{\widetilde{X}}, \; \; \widetilde{\mu}(\widetilde{\mu}^{-1}(a))\in B_{\widetilde{Y}}, \; \; \widetilde{\nu}(\widetilde{\mu}^{-1}(a)) = (b,c)(a,d)(c,d)[\widetilde{\mu}(\widetilde{\mu}^{-1}(a))] = d\in A_{\widetilde{Y}}, \text{ and }\\
    &\widetilde{\mu}^{-1}(d)\in A_{\widetilde{X}}, \; \; \widetilde{\mu}(\widetilde{\mu}^{-1}(d))\in A_{\widetilde{Y}}, \; \; \widetilde{\nu}(\widetilde{\mu}^{-1}(d)) = (b,c)(a,d)(c,d)[\widetilde{\mu}(\widetilde{\mu}^{-1}(d))] = b\in A_{\widetilde{Y}}.
    \end{split}
\end{equation*}

This shows that $\widetilde{\nu}(A_{\widetilde{X}})\cap A_{\widetilde{Y}} = (\widetilde{\mu}(A_{\widetilde{X}})\cap A_{\widetilde{Y}})\cup\{a\},$ so $|\widetilde{\nu}(A_{\widetilde{X}})\cap A_{\widetilde{Y}}|=|\widetilde{\mu}(A_{\widetilde{X}})\cap A_{\widetilde{Y}}|+1.$ This is a contradiction with the maximilaty of $|\widetilde{\mu}(A_{\widetilde{X}})\cap A_{\widetilde{Y}}|.$\\

As we reached a contradiction in both cases, it follows that in fact
$s = r-1.$ Therefore, there exists a sequence $\Sigma$ of $(X,Y)$-friendly swaps which does not involve $u$ or $v$ and transforms $\sigma$ into a bijection $\mu$ satisfying $\mu(A_X) = A_Y$ and $\mu(B_X) = B_Y.$ 
\end{proof}\\

\subsection{Endgame}
Let $\Sigma$ and $\mu$ be as in the statement of Lemma 6.2. Note that as $\Sigma$ does not involve $u$ and $v,$ it is true that $\mu^{-1}(u) = u'$ and $\mu^{-1}(v) = v'.$
Let $w$ be an arbitrary neighbor of $u$ in $B_Y.$ Denote $w' = \mu^{-1}(w).$ \\

From the assumption $\delta(X)+\delta(Y)\ge 3r/2+1,$ we conclude that the inequality
$$
|N(v)| + |\mu(N(v'))| +|N(w)|+|\mu(N(w'))|\ge 3r+2 
$$ holds
since $|N(v)|\ge \delta(Y), |N(w)|\ge \delta(Y)$ and 
$|\mu(N(v'))|\ge \delta(X),|\mu(N(w'))|\ge \delta(X).$ At the same time, all four sets $N(v), N(w),\mu(N(v'))$ and $\mu(N(w'))$ belong to $A_Y,$ which has cardinality $r.$ By the pigeonhole principle, we deduce that there exist at least two vertices $z_1$ and $z_2$ in $A_Y,$ which belong to all four sets $N(v), N(w),\mu(N(v'))$ and $\mu(N(w')).$ In other words, $z_1$ and $z_2$ satisfy: 
\vspace*{-0.4em}
\begin{itemize}
    \setlength\itemsep{-0.4em}
    \item $(z_1, v), (z_2,v) \in E(Y),$
    \item 
    $(\mu^{-1}(z_1), v'), (\mu^{-1}(z_2),v') \in E(X),$
    \item $(z_1, w), (z_2,w) \in E(Y),$ 
    \item $(\mu^{-1}(z_1), w'), (\mu^{-1}(z_2),w') \in E(X).$
\end{itemize}

We finish the proof by considering two cases.\\

\noindent
\textbf{Case I:} Suppose that one of $z_1$ and $z_2$ is $u.$ Then $u$ and $v$ are adjacent in $Y,$ so the $(X,Y)$-friendly swap $uv$ exchanges them.\\ 
\textbf{Case II:} Suppose that $z_1$ and $z_2$ are both distinct from $u.$ In that case, we can check that the sequence $\widetilde{\Sigma}$ of $(X,Y)$-friendly swaps
$$
\widetilde{\Sigma} = vz_1, wz_2, wz_1, wu,  wz_2, wz_1, vz_1, vz_2, wz_2
$$
exchanges $u$ and $v$ from $\mu.$ Therefore, the sequence $\Sigma^* = \Sigma,\widetilde{\Sigma},\rev(\Sigma)$ exchanges $u$ and $v$ from $\sigma.$

\begin{figure}[!htb]
\begin{center}
   \includegraphics[width = .8\linewidth]{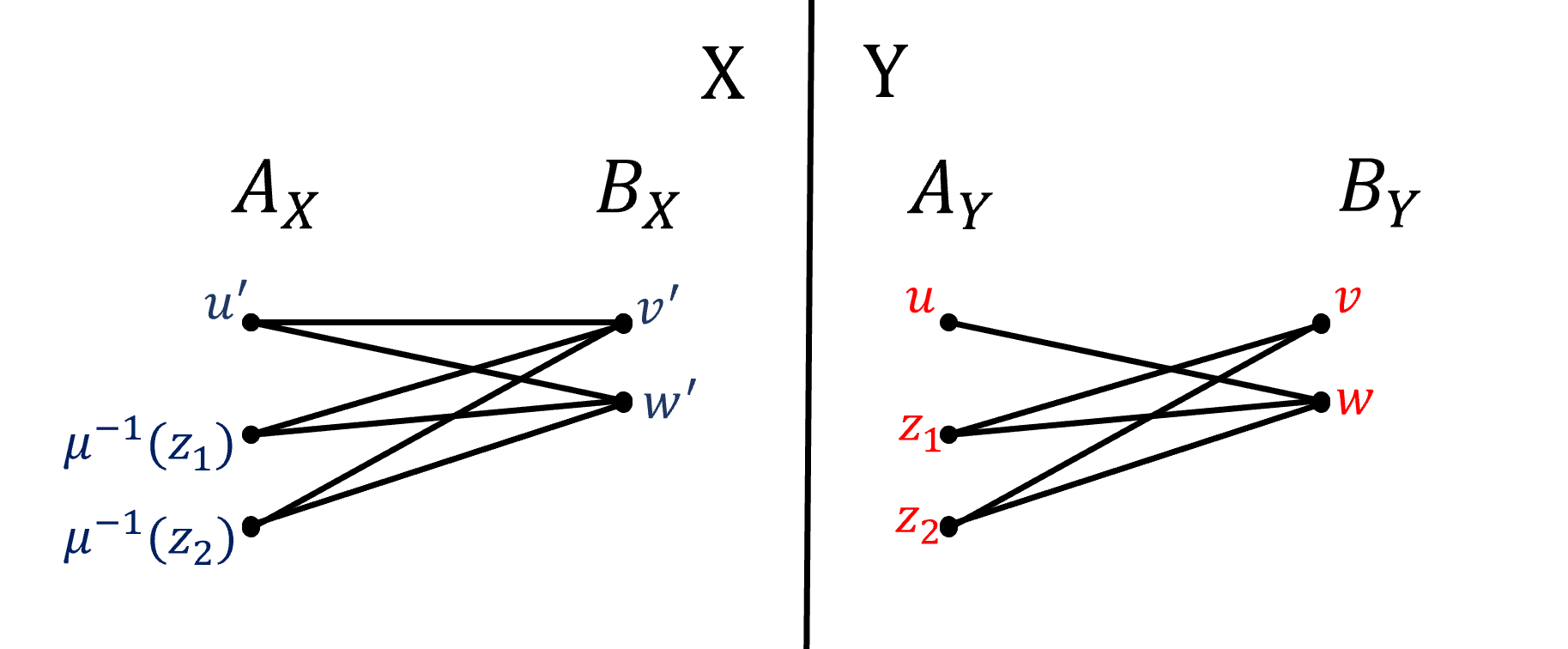}\caption{
   Relevant edges and vertices to Case II in Section 6.3
   }\label{fig:end_proof_k_rr_full}
\end{center}
\end{figure}
\section{Lower Bounds for Bipartite Graphs: Theorem 1.11}
\noindent
\textbf{Theorem 1.11}: \textit{ Let $r\ge 2$ and $\delta_1, \delta_2$ be non-negative integers satisfying
$$\delta_1 + \delta_2 =  \Big{\lfloor} \frac{3r}{2}\Big{\rfloor} ,\; \delta_1 \le r,\;  \delta_2\le r.$$ 
Then there exist two edge-subgraphs $X$ and $Y$ of $K_{r,r}$ such that $\delta(X) = \delta_1, \delta(Y) = \delta_2,$ and $\FS(X,Y)$ has more than two connected components.}\\
\begin{proof} First, note that the above inequalities imply $\lfloor r/2\rfloor \le \delta_1\le r$ and $\lfloor r/2\rfloor \le \delta_2\le r.$\\

Let the partite sets in the bipartion of $X$ be $A_X\sqcup B_X$ and $C_X\sqcup D_X,$ where $|A_X| =|D_X| =\lceil r/2\rceil $ and $|B_X| = |C_X|=\lfloor r/2\rfloor.$
Similarly, let the partite sets in the bipartion of $Y$ be $A_Y\sqcup C_Y$ and $B_Y\sqcup D_Y,$ where $|A_Y| =|D_Y| =\lceil r/2\rceil $ and $|B_Y| = |C_Y|=\lfloor r/2\rfloor.$ Fix a bijection $\sigma : V(X)\longrightarrow V(Y)$ such that $\sigma(A_X) = A_Y, \sigma(B_X) = B_Y,\sigma(C_X) = C_Y,$ and $\sigma(D_X) = D_Y.$\\

We first construct $X.$ In $X,$ every vertex in $A_X$ is adjacent to every vertex in $C_X$ and every vertex in $B_X$ is adjacent to every vertex in $D_X.$ Every vertex in $A_X$ is adjacent to $\delta_1 - |C_X| = \delta_1 - \lfloor r/2 \rfloor$ vertices in $D_X$ and every vertex in $D_X$ is adjacent to $\delta_1 - |B_X| = \delta_1-\lfloor r/2 \rfloor$ vertices in $A_X.$ This is easy to achieve as $|A_X| = |D_X| = \lceil r/2 \rceil$ and $0 \le \delta_1 - \lfloor r/2\rfloor\le \lceil r/2\rceil.$ Similarly, every vertex in $B_X$ is adjacent to $\max(0,\delta_1 - \lceil r/2\rceil)$ vertices in $C_X$ and every vertex in $C_X$ is adjacent to $\max(0,\delta_1 - \lceil r/2\rceil)$ vertices in $B_X.$ Clearly, $X$ is a bipartite graph with minimum degree at least $\delta_1.$\\

Now, we construct $Y$ in much the same way. In $Y,$ every vertex in $A_Y$ is adjacent to every vertex in $B_Y$ and every vertex in $C_Y$ is adjacent to every vertex in $D_Y.$ For any two vertices $u\in A_Y$ and $v\in D_Y,$ we connect $u$ and $v$ if and only if $(\sigma^{-1}(u), \sigma^{-1}(v))$ is not an edge in $X.$ Using that $\sigma(A_X) = A_Y,\sigma(D_X) = D_Y$ and the construction of $X,$ every vertex of $A_Y$ is adjacent to $$\lceil r/2\rceil - (\delta_1-\lfloor r/2\rfloor) = r - \delta_1 = \delta_2 - \lfloor r/2\rfloor$$ vertices in $D_Y$ and every vertex of $D_Y$ is adjacent to $\delta_2 - \lfloor r/2\rfloor$ vertices in $A_Y.$ Similarly, for any two vertices $u\in B_Y,v\in C_Y,$ we connect $u$ and $v$ if and only if $(\sigma^{-1}(u), \sigma^{-1}(v))$ is not an edge in $X.$ Thus, each vertex in $B_Y$ has 
$$
|C_Y| - \max(0,\delta_1 - \lceil r/2\rceil) = 
\lfloor r/2\rfloor - \max(0,\delta_1 - \lceil r/2\rceil)
$$
neighbors in $C_Y.$ So, the degree of every vertex in $B_Y$ is
$$
|A_Y| + \lfloor r/2\rfloor - \max(0,\delta_1 - \lceil r/2\rceil) = 
r-\max(0,\delta_1 - \lceil r/2\rceil).
$$
We claim that this is always at least $\delta_2.$ Indeed, $r-0\ge \delta_2$ and $r - (\delta_1-\lceil r/2\rceil) = \lceil 3r/2\rceil-\delta_1\ge \delta_2.$ The same holds for $C_Y.$
The constructed graph $Y$ is bipartite and has minimum degree at least $\delta_2.$\\

To show that $\FS(X,Y)$ has at least three connected components, note that $\sigma$ is an isolated vertex because for no $u,v\in V(Y)$ is it the case that 
$(u,v)\in E(Y)$ and $(\sigma^{-1}(u), \sigma^{-1}(v))\in E(X)$ simultaneously hold. Now, partition the bijections from $V(X)$ to $V(Y)$ as follows:
$$
K_0 = \{\sigma\},
$$
$$
K_1 = \{\tau\neq \sigma : sign(\tau^{-1}\circ\sigma) \equiv
|\sigma(A_X\cup B_X)\cap (A_Y\cup C_Y)| - |\tau(A_X\cup B_X)\cap (A_Y\cup C_Y)|  \pmod{2}\},
$$
$$
K_2 = \{\tau\neq \sigma : sign(\tau^{-1}\circ\sigma) \not\equiv
|\sigma(A_X\cup B_X)\cap (A_Y\cup C_Y)| - |\tau(A_X\cup B_X)\cap (A_Y\cup C_Y)|  \pmod{2}\}.
$$
As $r\ge 2,$ it is easy to show that none of the three sets is empty. However, there are no edges from a vertex in $K_0$ to a vertex in $K_1$ or $K_2$ since $\sigma$ is isolated. There are no edges between a vertex in $K_1$ and a vertex in $K_2$ by Proposition 2.2. Thus, $\FS(X,Y)$ has at least three connected components.\end{proof}\\

\section{Future Work}
One future direction is finding tight estimates in Problem 1.3. Likely, this would require both improved lower and upper bounds. Based on the results in the current paper, we make the following conjectures.\\

\begin{figure}[!htb]
\begin{center}
   \includegraphics[width = 
   \linewidth]{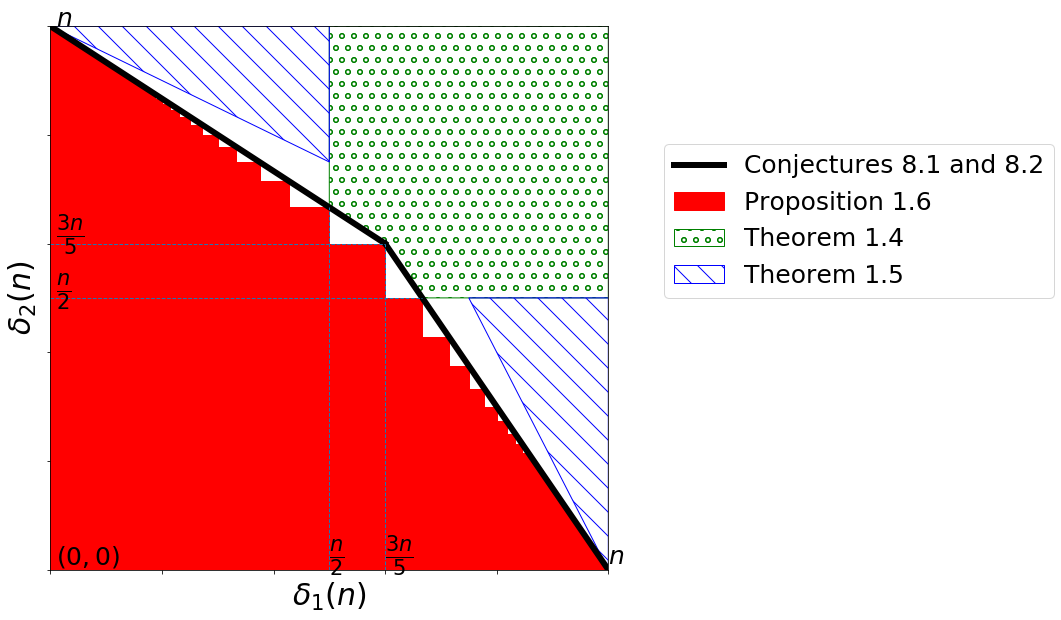}\caption{Visualization of proved and conjectured estimates regarding Problem 1.3 in the current paper. In solid red are pairs $(\delta_1(n),\delta_2(n))$ which do not guarantee that $\FS(X,Y)$ is connected. In dotted green and striped blue are pairs which guarantee that $\FS(X,Y)$ is connected. White portions are still unresolved.
   The black line corresponds to the conjectured tight bound. Note that the ``corners'' of the shape given by Proposition 1.6 lie on that line, as does part of the boundary of the shape given by Theorem 1.4.
   \textbf{Remark: }Additive constants are omitted in the diagram.}\label{fig:conj_diagram}
\end{center}
\end{figure}

\noindent
\textbf{Conjecture 8.1}: \textit{Suppose that $X$ and $Y$ are two connected $n$-vertex graphs satisfying the condition
    $$2\min(\delta(X),\delta(Y))+3\max(\delta(X), \delta(Y))\ge 3n.$$
Then $\FS(X,Y)$ is connected.}\\

We further conjecture that this is tight up to an additive constant:\\

\noindent
\textbf{Conjecture 8.2}: \textit{There exists a constant $c>0$ such that for any three positive integers $n,\delta_1, \delta_2$ satisfying
$$
2\min(\delta_1,\delta_2)+3\max(\delta_1, \delta_2)\le 3n-c,\; \; \delta_1\le n-1, \; \; \delta_2\le n-1,
$$
there exist two connected simple $n$-vertex graphs $X$ and $Y$ with $\delta(X) = \delta_1$ and $\delta(Y)=\delta_2$ such that $\FS(X,Y)$ is disconnected.}\\

See \autoref{fig:conj_diagram} for an illustration of proven and conjectured estimates regarding Problem 1.3.\\

Finally, we propose the following natural hypergraph generalization of friends-and-strangers\\ graphs. As usual, in a $k$-uniform hypergraph $G,$ each edge $e\in E(G)$ is a set of $k$ distinct vertices of $G.$ Using this notation, we make the following definition.\\

\newpage

\noindent
\textbf{Definition 8.3}: For two $k$-uniform hypergraphs $X$ and $Y$ with the same (finite) number of vertices, their friends-and-strangers graph $\FS_k(X,Y)$ is defined as follows:
\vspace*{-0.8em}
\begin{itemize}
    \setlength\itemsep{-0.4em}
    \item The vertex set $V(\FS_k(X,Y))$ of $\FS_k(X,Y)$ consists of all bijections $\sigma: V(X)\longrightarrow V(Y).$
    \item Two distinct bijections $\sigma, \tau \in V(\FS(X,Y))$ are adjacent if and only if there exist two edges $e\in E(X)$ and $f\in E(Y)$ such that:
    \vspace*{-0.8em}
    \begin{itemize}
    \setlength\itemsep{-0.4em}
        \item $\sigma(e) = f,$
        \item There exists a permutation $\mu: V(Y)\longrightarrow V(Y),$ equal to the identity on $V(Y)\backslash f,$ for which $\sigma = \mu \circ \tau.$
    \end{itemize}
\end{itemize}

When $k=2$ this reduces to Definition 1.1. We pose the following analogue of Problem 1.3:\\

\noindent
\textbf{Problem 8.4}: \textit{Find conditions on the pairs $(\delta^k_1(n),\delta^k_2(n))$ which guarantee that if $X$ and $Y$ are any two $k$-uniform hypergraphs such that each vertex in $V(X)$ is in at least $\delta^k_1(n)$ edges of $X$ and each vertex in $V(Y)$ is in at least $\delta^k_2(n)$ edges of $Y,$ then
$\FS_k(X,Y)$ is connected.}\\

One could also study the connectivity of $\FS_k(X,Y)$ when one of the hypergraphs has a specific fixed structure or $X$ and $Y$ are randomly generated. Similarly, one can try to identify structural obstructions to the connectivity of $\FS_k(X,Y)$ such as the one given in Proposition 2.2 for the special case $k=2.$

\section*{Acknowledgements}
This research was conducted at the University of Minnesota Duluth Mathematics REU and was supported, in part,  by NSF-DMS Grant 1949884 and NSA Grant H98230-20-1-0009. I was fully supported by the Department of Mathematics at Princeton University.
I want to thank Joe Gallian for organizing the REU program and introducing me to the problem. I am also grateful to Colin Defant and Noah Kravitz for  many helpful discussions on friends-and-strangers graphs. Finally, I want to express my gratitude towards
Colin Defant, Joe Gallian, Noah Kravitz, Mitchell Lee, and Nina Zubrilina for reading this paper at different stages of its development and suggesting many invaluable improvements.

\bibliography{main}{}
\bibliographystyle{elsarticle-num}

\vspace*{1cm}
 \textsc{Department of Mathematics, Princeton University,
    Princeton, New Jersey, 08544}\par\nopagebreak
  \textit{Email address:} \href{mailto:kirilb@princeton.edu}{kirilb@princeton.edu} 
\end{document}